\documentclass[12pt]{article}
\usepackage{amssymb,amsfonts,amsmath,psfrag,eepic,colordvi,graphicx,epsfig}
\usepackage{amssymb,latexsym,graphics,array}
\usepackage{MnSymbol}
\usepackage[enableskew]{youngtab}
\usepackage{float,enumerate,enumitem}
\usepackage{tikz,ifthen}
\usepackage{pgfplots}
\usetikzlibrary{math}

\usepackage{framed}
\usepackage{hyperref}
\hypersetup{colorlinks=true}

\topskip  
\numberwithin{equation}{section}
\newtheorem{theorem}{Theorem}[section]
\newtheorem{proposition}[theorem]{Proposition}

\newtheorem{definition}[theorem]{Definition}

\newtheorem{lemma}[theorem]{Lemma}

\newtheorem{example}[theorem]{Example}

\newtheorem{observation}[theorem]{Observation}

\begin{document}
	\parskip 6pt
	
	\pagenumbering{arabic}
	\def\sof{\hfill\rule{2mm}{2mm}}
		\def\des{\mathsf{des}}
		\def\Si{\mathsf{Si}}
		\def\Oint{\mathsf{Oint}}
		\def\Yleaf{\mathsf{Yleaf}}
	 		\def\lleaf{\mathsf{lleaf}}
		\def\rleaf{\mathsf{rleaf}}
		\def\Rleaf{\mathsf{Rleaf}}
	 		\def\oleaf{\mathsf{oleaf}}
	 			\def\Oleaf{\mathsf{Oleaf}}
		 		\def\yleaf{\mathsf{yleaf}}
		\def\oint{\mathsf{oint}}
		\def\lint{\mathsf{lint}}
		\def\si{\mathsf{si}}
			\def\rpd{\mathsf{rpd}}
				\def\pd{\mathsf{pd}}
			\def\sp{\mathsf{sp}}
			\def\lap{\mathsf{lap}}
	\def\asc{\mathsf{asc}}
	\def\exc{\mathsf{exc}}
	\def\cyc{\mathsf{cyc}}
	\def\ap{\mathsf{ap}}
	\def\ls{\leq}
	\def\gs{\geq}
	\def\SS{\mathcal S}
	\def\qq{{\bold q}}
	\def\MM{\mathcal M}
	\def\TT{\mathcal T}
	\def\EE{\mathcal E}
	\def\lsp{\mbox{lsp}}
	\def\rsp{\mbox{rsp}}
	\def\pf{\noindent {\it Proof.} }
	\def\mp{\mbox{pyramid}}
	\def\mb{\mbox{block}}
	\def\mc{\mbox{cross}}
	\def\qed{\hfill \rule{4pt}{7pt}}
	\def\block{\hfill \rule{5pt}{5pt}}
	\def\lr#1{\multicolumn{1}{|@{\hspace{.6ex}}c@{\hspace{.6ex}}|}{\raisebox{-.3ex}{$#1$}}}
	\def\red{\textcolor{red}}
	
 \allowdisplaybreaks 

\begin{center}
{\Large\bf Combinatorics on bi-$\gamma$-positivity of $1/k$-Eulerian polynomials}
\end{center}

\begin{center}
	
	{\small Sherry H.F. Yan$^a$, Xubo Yang$^a$,   Zhicong Lin$^{b,*}$\footnote{$^*$Corresponding author. E-mail address: linz@sdu.edu.cn} }
	
	$^{a}$Department of Mathematics,
	Zhejiang Normal University\\
	Jinhua 321004, P.R. China
	
	$^{b}$ Research Center for Mathematics and Interdisciplinary Sciences,  Shandong University \& Frontiers Science Center for Nonlinear Expectations, Ministry of Education, Qingdao 266237, P.R. China

\end{center}

\noindent {\bf Abstract.}
 The $1/k$-Eulerian polynomials $A^{(k)}_{n}(x)$  were introduced as ascent   polynomials over $k$-inversion sequences   by Savage  and  Viswanathan. 
The bi-$\gamma$-positivity of the  $1/k$-Eulerian polynomials $A^{(k)}_{n}(x)$ was known but to give a 
  combinatorial interpretation of the corresponding    bi-$\gamma$-coefficients   still remains open. The study of the theme of bi-$\gamma$-positivity   from purely  combinatorial aspect was proposed by Athanasiadis. In this paper, we provide a  combinatorial interpretation for the  bi-$\gamma$-coefficients of $A^{(k)}_{n}(x)$ by using the model of certain ordered labeled  forests. Our combinatorial approach consists of three main steps: 
  \begin{itemize}
  	\item    construct a bijection between $k$-Stirling permutations and certain forests that are named   increasing pruned even $k$-ary forests; 
  \item  introduce a generalized Foata–Strehl  action on increasing pruned even $k$-ary trees which implies the  longest ascent-plateau polynomials  over  $k$-Stirling permutations with initial letter $1$ are $\gamma$-positive, a result that may have independent interest; 
  \item develop two crucial transformations on increasing pruned even $k$-ary forests to conclude our combinatorial interpretation. 
  \end{itemize}




\section{Introduction}

Following Savage  and  Viswanathan \cite{Savage-2012},   the {\em$1/k$-Eulerian polynomials} $A^{(k)}_n(x)$ are defined by 
$$
\sum_{n\geq 0}A_n^{(k)}(x){z^n\over n!}=\left({1-x\over e^{kz(x-1)}-x}  \right)^{{1\over k}}.
$$
When $k=1$,  the $1/ k$-Eulerian polynomials $A^{(k)}_n(x)$  reduce to the classical Eulerian polynomials 
$$
A_n(x)=\sum_{\pi\in \mathfrak{S}_n}x^{\des(\pi)}=\sum_{\pi\in \mathfrak{S}_n}x^{\exc(\pi)},
$$ 
where   $\mathfrak{S}_n$ denotes the set of all permutations of $[n]:=\{1, 2, \ldots, n\}$ and $\des(\pi)$ (resp., $\exc(\pi)$) denotes the number of descents (resp., excedances) of $\pi$. 
Savage  and  Viswanathan~\cite{Savage-2012} discovered the combinatorial interpretation of $A_n^{(k)}(x)$ in terms of ascent polynomials over $k$-inversion sequences
$\{(e_1, e_2, \ldots, e_n)\in \mathbb{Z}^n\mid 0\leq  e_i\leq (i-1)k \}$, which are closely related to the generalized lecture hall partitions~\cite{Savage-2016}. They also 
 proved that 
 \begin{equation}
  A^{(k)}_n(x)= \sum_{\pi\in \mathfrak{S}_n}x^{\exc(\pi)}k^{n-\cyc(\pi)},
 \end{equation}
where $\cyc(\pi)$ denotes the number of cycles  in the disjoint cycle representation of $\pi$. 
A bijection linking the interpretation of $A^{(k)}_n(x)$ as $k$-inversion sequences with that as permutations was constructed by Chao, Ma, Ma and Yeh~\cite{JM2}.

 Let $M=\{1^{k_1}, 2^{k_2}, \ldots, n^{k_n}\}$ be a multiset,  where $i^{k_i}$ signifies  $k_i$ occurrences of $i$. 
A permutation $\pi$ of  $M$ is said to be  a {\em Stirling   permutation }  if  $\pi_i=\pi_k$ implies that  $\pi_j\geq \pi_i$ for all $i<j<k$.
For $k\geq 1$, 
given a set $A=\{a_1, a_2, \ldots, a_m\}$ with $a_1<a_2<\ldots<a_m$, let  $A_k$ denote
the multiset $\{a_1^k, a_2^k, \ldots, a_n^k\}$.
Stirling permutations were originally introduced by Gessel and Stanley \cite{Gessel-Stanley-1978} in the case of the multiset $[n]_k$. 
A Stirling  permutation $\pi$ on  $[n]_k$ is said to be a {\em $k$-Stirling permutation }of order $n$. Let $\mathcal{Q}_{n}(k)$ denote the set of $k$-Stirling permutations  of order $n$.
Given $\pi=\pi_1\pi_2\cdots \pi_{kn}\in \mathcal{Q}_{n}(k)$, an index $i$    is said to be a {\em longest   plateau} if
$$
\pi_i=\pi_{i+1}=\cdots=\pi_{i+k-1},
$$
and 
 an index $i$  is said to be a {\em longest  ascent-plateau} if
$$
\pi_i<\pi_{i+1}=\pi_{i+2}=\cdots=\pi_{i+k}.
$$ 
A {\em left longest   ascent-plateau}  of $\pi$ is a longest ascent-plateau of $\pi$ patched with a zero in the front of $\pi$.  
Let $\ap(\pi)$ and
$
\lap(\pi)
$  denote the number  of  longest ascent-plateaux and the number of left   longest  ascent-plateaux  of $\pi$, respectively.  
In \cite{Ma-Mansour-2015}, Ma and Mansour proved that
\begin{equation}
	 A^{(k)}_n(x)=\sum_{\pi\in \mathcal{Q}_n(k)}x^{\ap(\pi)}, \,\,\,\,  x^nA_n^{(k)}({1/x})= \sum_{\pi\in \mathcal{Q}_n(k)}x^{\lap(\pi)}.
\end{equation}
Through bijections with $k$-Stirling permutations, Liu~\cite{Liu} found certain statistics on $k$-inversion sequences and trapezoidal words that also give rise to the $1/ k$-Eulerian polynomials $A^{(k)}_n(x)$. 

Let $h(x)=\sum_{i=0}^{n}a_ix^i$ be a polynomial with nonnegative real coefficients. 
We say that $h(x)$ is {\em unimodal } if
$$a_0\leq a_1\leq a_2\leq\cdots \leq a_m\geq a_{m+1}\geq \cdots\geq a_n $$
for some $m$,  where the index $m$ is called the {\em mode} of $h(x)$. 
The polynomial $h(x)$ is said  to be  {\em alternating increasing} if $$a_0\leq a_n\leq a_1\leq a_{n-1}\leq \cdots \leq a_{\lfloor {n+1\over 2}\rfloor}. $$ Clearly, alternating increasing is a stronger property than unimodality. 
If $a_i=a_{n-i}$ for all $0\leq i\leq n$, then  $h(x)$ is said to be {\em symmetric}. 
 We say that $h(x)$ is {\em $\gamma$-positive} if it can be expanded as
$$
h(x)=\sum_{i=0}^{\lfloor {n\over 2}\rfloor} \gamma_k x^k(1+x)^{n-2k}
$$
with the $\gamma$-coefficient $\gamma_k\geq 0$.  Clearly,   $\gamma$-postivity    implies symmetry and unimodality. 
Foata and Sch\"uzenberger \cite{FS-70} proved that the Eulerian polynomial $A_n(x)$ is   $\gamma$-positive and provided a combinatorial interpretation of the corresponding $\gamma$-coefficients in terms of permuations without double descents.     An elegant combinatorial proof   via a group action was later constructed by Foata and Strehl \cite{Foata-Strehl-1974}. Since then, $\gamma$-positive polynomials arising in combinatorics and geometry have been extensively investigated; see \cite{A-2018, Branden-2008, br2, Chen-2023, Chow-2008, Fu, Ji-Lin, Lin, Lin-Ma-Zhang-2021,  Lin-Zeng, JM, Ma-Ma-Yeh-2019, Ma-Ma-Yeh-2024, Yan-Huang-Yang-2021} and the references therein.

\begin{proposition}\label{prop1} {\upshape(See \cite{Branden-2021})}
	Let $h(x)$ be a polynomial of degree $n$. Then $h(x)$ can be uniquely decomposed as $h(x)=a(x)+xb(x)$, where $a(x)$ is symmetric with degree $n$ and $b(x)$ is symmetric  with degree less than $n$. More precisely, we have
	$$
	a(x)={h(x)-x^{n+1}h(1/x)\over 1-x}\,\,\,\,\,\, \mbox{and}\,\,\,\, b(x)={x^{n}h(1/x)-h(x)\over 1-x}.	$$
\end{proposition}
The decomposition  $h(x)=a(x)+xb(x)$ in Proposition \ref{prop1} is called the {\em symmetric decomposition} of $h(x)$.  As observed in \cite{Branden-2021},  the polynomial $h(x)$ is alternating increasing if and only if $a(x)$ and $b(x)$ are both unimodal and have nonnegative coefficients. If $a(x)$ and $b(x)$ are both $\gamma$-positive, then we say that $h(x)$ is {\em bi-$\gamma$-positive}. Note that if $h(x)$ is bi-$\gamma$-positive, then it is alternating increasing. Bi-$\gamma$-positivity  naturally extends 
$\gamma$-positivity from symmetric polynomials to general polynomials. There has been recent interest in studying bi-$\gamma$-positive polynomials  arising in combinatorics and geometry;  see \cite{A-2018,A-2020, Branden-2021, DHL, Han, Lin-Xu-Zhao} and the references therein. The study of the theme of bi-$\gamma$-positivity   from purely  combinatorial aspect was proposed by Athanasiadis~\cite{A-2018}.

Define 
$$
\overline{\mathcal{Q}}_{n}(k)=\{\pi=\pi_1\pi_2\ldots \pi_{kn}\in \mathcal{Q}_n(k)\mid \pi_1=\pi_2=\ldots =\pi_k \}
$$
and 
$$
\widehat{\mathcal{Q}}_{n}(k)= \mathcal{Q}_{n}(k)-\overline{\mathcal{Q}}_{n}(k).
$$
 In~\cite{Ma-Ma-Yeh-2024}, Ma-Ma-Yeh-Yeh provided the following combinatorial interpretation of the symmetric decomposition of $A^{(k)}_{n}(x)$. 
 \begin{proposition}{\upshape(See~\cite[Proposition 4.2]{Ma-Ma-Yeh-2024})}\label{prop2} 
 Suppose that the symmetric decomposition of $A^{(k)}_{n}(x)$ is given by $a^{(k)}_{n}(x)+xb^{(k)}_{n}(x)$, then we have 
 $$
 a^{(k)}_{n}(x)=\sum_{\pi\in \overline{\mathcal{Q}}_n(k) }x^{\ap(\pi)}\,\,\,\,\, \mbox{and}\,\,\,\, xb^{(k)}_{n}(x)=\sum_{\pi\in \widehat{\mathcal{Q}}_{n}(k) }x^{\ap(\pi)}.
 $$ 
 \end{proposition}

Via the machine of context-free grammars,  the  $1/k$-Eulerian polynomials $A^{(k)}_{n}(x)$ were also proved to be bi-$\gamma$-positive~\cite{Ma-Ma-Yeh-2024}. Moreover, they provided a combinatorial interpretation for  the bi-$\gamma$-coefficients  of   $A^{(2)}_{n}(x)$ in terms of permutations.  The combinatorial interpretation for the    bi-$\gamma$-coefficients  of    $A^{(k)}_{n}(x)$ for general $k$ still remains open.

For $\gamma$-positive or bi-$\gamma$-positive polynomials arising in combinatorics or geometry, an interesting and challenging problem is to find combinatorial interpretations for their  corresponding  $\gamma$-coefficients or bi-$\gamma$-coefficients. 
 Ordered labeled trees play an important role in the  combinatorial interpretation of the  $\gamma$-coefficients or   bi-$\gamma$-coefficients of enumerative polynomials; see  for instance~\cite{Chen-2023, Lin, Lin-Xu-Zhao, Yan-Huang-Yang-2021}. The main objective of this paper is to provide  a combinatorial interpretation for the    bi-$\gamma$-coefficients  of    $A^{(k)}_{n}(x)$ by using the model of certain   ordered labeled   forests.

\section{Main result}
  In this section, we shall state our main result whose proof will be  presented in the subsequent sections.

Let us first review some terminology related to trees.
A {\em tree} is an acyclic connected graph, and a {\em forest} is a graph such that every connected component is a tree. 
An  {\em ordered} tree is a tree with one designated node, which is called the root, and the subtrees of
each node are linearly ordered. 
In this paper we will assume that all the trees are ordered.  In a  tree $T$,   the {\em level} of a node $v$ in $T$  is defined to be the length of the unique path from the root to $v$.  For each node of the tree $T$, we say that it is of degree $k$
  if it has exactly $k$ 
  children.  The   nodes of degree $0$ of the tree $T$ are called {\em leaves} and  the nodes that are not leaves are called {\em internal nodes}.  A  tree $T$ is said to be an  {\em even  $k$-ary tree} if  all nodes on even levels  are of degree $k$. Given an even $k$-ary tree $T$,  for each node $u$ on even level, if  all the children of $u$ are leaves, we will chop off all the leaves  attached to $u$.    
    Intuitively,   the resulting tree is called    a {\em pruned even $k$-ary  tree}. 
    
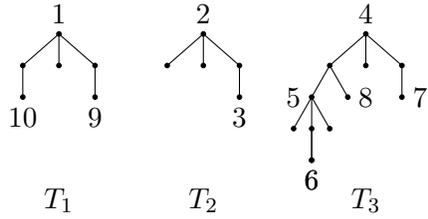
\begin{figure}
	\begin{center}
		\begin{tikzpicture}[font =\small , scale = 0.3]
			\filldraw[fill=black](5*0.8,0.7*0)circle(0.1);
			\coordinate [label=above:$1$] (x) at (5*0.8,0.7*0);
			\filldraw[fill=black](3*0.8,0.7*-2)circle(0.1);
			
			\filldraw[fill=black](3*0.8,0.7*-4)circle(0.1);
			\coordinate [label=below:$10$] (x) at (3*0.8,0.7*-4);
			
			\draw(3*0.8,0.7*-2)--(3*0.8,0.7*-4);
			\filldraw[fill=black](5*0.8,0.7*-2)circle(0.1);
			
			\filldraw[fill=black](7*0.8,0.7*-2)circle(0.1);
			
			\filldraw[fill=black](7*0.8,0.7*-4)circle(0.1);
			\coordinate [label=below:$9$] (x) at (7*0.8,0.7*-4);
			\draw(7*0.8,0.7*-2)--(7*0.8,0.7*-4);
			\draw(5*0.8,0.7*0)--(3*0.8,0.7*-2);
			\draw(5*0.8,0.7*0)--(5*0.8,0.7*-2);
			\draw(5*0.8,0.7*0)--(7*0.8,0.7*-2);
			
			\coordinate [label=above:$ T_1$] (x) at (5*0.8,0.7*-12);
			
			\filldraw[fill=black](13*0.8,0.7*0)circle(0.1);
			\coordinate [label=above:$2$] (x) at (13*0.8,0.7*0);
			\filldraw[fill=black](13*0.8,0.7*-2)circle(0.1);
			
			\filldraw[fill=black](11*0.8,0.7*-2)circle(0.1);
			
			\filldraw[fill=black](15*0.8,0.7*-2)circle(0.1);
			\filldraw[fill=black](15*0.8,0.7*-4)circle(0.1);
			\coordinate [label=below:$3$] (x) at (15*0.8,0.7*-4);
			\draw(13*0.8,0.7*0)--(13*0.8,0.7*-2);
			\draw(13*0.8,0.7*0)--(11*0.8,0.7*-2);
			\draw(13*0.8,0.7*0)--(15*0.8,0.7*-2);
			\draw(15*0.8,0.7*-4)--(15*0.8,0.7*-2);
			
			\coordinate [label=above:$ T_2$] (x) at (13*0.8,0.7*-12);
			
			\filldraw[fill=black](22*0.8,0.7*0)circle(0.1);
			\coordinate [label=above:$4$] (x) at (22*0.8,0.7*0);
			\filldraw[fill=black](20*0.8,0.7*-2)circle(0.1);
			
			\filldraw[fill=black](19*0.8,0.7*-4)circle(0.1);
			\coordinate [label=left:$5$] (x) at (19*0.8,0.7*-4);
			\filldraw[fill=black](21*0.8,0.7*-4)circle(0.1);
			\coordinate [label=right:$8$] (x) at (21*0.8,0.7*-4);
			\draw(19*0.8,0.7*-4)--(20*0.8,0.7*-2);
			\draw(21*0.8,0.7*-4)--(20*0.8,0.7*-2);
			\filldraw[fill=black](19*0.8,0.7*-6)circle(0.1);
			
			\filldraw[fill=black](18*0.8,0.7*-6)circle(0.1);
			
			\filldraw[fill=black](20*0.8,0.7*-6)circle(0.1);
			
			\draw(19*0.8,0.7*-4)--(18*0.8,0.7*-6);
			\draw(19*0.8,0.7*-4)--(19*0.8,0.7*-6);
			\draw(19*0.8,0.7*-4)--(20*0.8,0.7*-6);
			\filldraw[fill=black](19*0.8,0.7*-8)circle(0.1);
			\coordinate [label=below:$6$] (x) at (19*0.8,0.7*-8);
			\draw(19*0.8,0.7*-6)--(19*0.8,0.7*-8);
			\filldraw[fill=black](22*0.8,0.7*-2)circle(0.1);
			
			\draw(22*0.8,0.7*0)--(22*0.8,0.7*-2);
			\draw(22*0.8,0.7*0)--(20*0.8,0.7*-2);
			\draw(22*0.8,0.7*0)--(24*0.8,0.7*-2);
			\filldraw[fill=black](24*0.8,0.7*-2)circle(0.1);
			
			\filldraw[fill=black](24*0.8,0.7*-4)circle(0.1);
			\coordinate [label=right:$7$] (x) at (24*0.8,0.7*-4);
			\draw(24*0.8,0.7*-2)--(24*0.8,0.7*-4);
			\coordinate [label=below:$6$] (x) at (19*0.8,0.7*-8);
			\filldraw[fill=black](19*0.8,0.7*-8)circle(0.1);			
			\draw(19*0.8,0.7*-6)--(19*0.8,0.7*-8);
			
			\coordinate [label=above:$ T_3$] (x) at (22*0.8,0.7*-12);
		\end{tikzpicture}
	\end{center}
	\caption{A forest  $F\in \mathcal{F}_{10}(3)$. \label{F}}
\end{figure}

 \begin{definition}[Increasing pruned  even $k$-ary tree]
 	Given a set $M$ with distinct positive integers, a labeled   tree $T$ is said to be an   increasing pruned  even $k$-ary tree  on $M$ if $T$ is a  pruned even $k$-ary tree in which the nodes on even levels are labeled by the elements of   $M$ satisfying that
 	\begin{itemize} 		 
 			\item  the labels of the nodes on the path from the root to any leaf is  increasing from top to bottom;
 		\item   the labels of the children of  each node  is   increasing from left to right.  		
 	\end{itemize}
 	\end{definition}
\begin{definition}[Increasing pruned  even $k$-ary forest]
For a set $M$, 	we say that a labeled forest $F$ is an   increasing pruned  even $k$-ary forest if each connected component is   an   increasing pruned even $k$-ary  tree, the labels of the roots are increasing from left to right and the labels of the  trees form a set partition of $M$. 
\end{definition}

Let $\mathcal{T}_{M}(k)$ and $\mathcal{F}_{M}(k)$ denote the set of    increasing pruned  even $k$-ary trees  and the set of  increasing pruned even $k$-ary forests on $M$, respectively.  When $M=[n]$, we simply write $\mathcal{T}_{M}(k)$ (resp.,~$\mathcal{F}_{M}(k)$) as $\mathcal{T}_{n}(k)$ (resp.,~$\mathcal{F}_{n}(k)$) for short.   In what follows, we will always write a  forest $F$ as a sequence  $(T_1, T_2, \ldots, T_m)$ for some $m$, where $T_i$ is the $i$-th tree of $F$ (counting from left to right).
  Fig.~\ref{F} displays  a forest  $F=(T_1, T_2, T_3)\in \mathcal{F}_{
	10 }(3)$.

In a tree $T\in   \mathcal{T}_{M}(k)$, a labeled node $u$  is called the {\em grand parent} of the labeled  node  $v$  if  one of a child of $u$ is the parent of $v$.   For instance, in the tree  $T_3$ of Fig.~\ref{F}, the nodes   $5,8$ and $7$ are grand children of  $4$.     A labeled   node $v$ (other than the root)  of the tree  $T$   is  {\em old } if  $v$ has the greatest label among all the grand children of its grand parent; otherwise, it is {\em young}.  The root of a tree is assumed to be   neither young nor old.  For example, in the tree  $T_3$ of Fig.~\ref{F}, the nodes $6$ and $8$ are old,  whereas the nodes $5$ and $7$  are young. 
Let $\lleaf(F)$  denote the total number of labeled leaves of the trees in $F$.  In a forest $F$, a tree $T$ of $F$ is called a {\em singleton} if   $T$ has only one node.    Denote by $\si(F)$ the number of singletons  of $F$.

  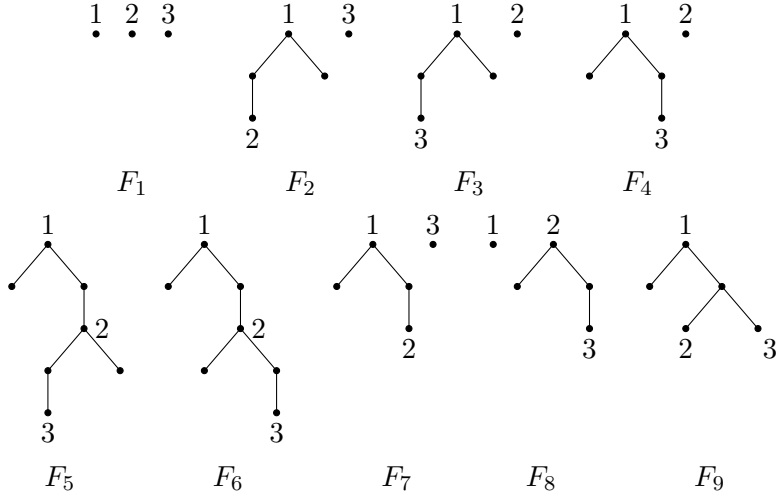
\begin{figure}
 	\begin{center}
 		\begin{tikzpicture}[font =\small , scale = 0.4]
 			\filldraw[fill=black](0.5*0.8,0.7*0)circle(0.1);
 			\filldraw[fill=black](2*0.8,0.7*0)circle(0.1);
 			\filldraw[fill=black](3.5*0.8,0.7*0)circle(0.1);
 			\coordinate [label={[text=black]above:$1$}] (x) at (0.5*0.8,0.7*0);
 			\coordinate [label={[text=black]above:$2$}] (x) at (2*0.8,0.7*0);
 			\coordinate [label={[text=black]above:$3$}] (x) at (3.5*0.8,0.7*0);
 			
 			\filldraw[fill=black](8.5*0.8,0.7*0)circle(0.1);
 			\filldraw[fill=black](7*0.8,0.7*-2)circle(0.1);
 			\filldraw[fill=black](10*0.8,0.7*-2)circle(0.1);
 			\filldraw[fill=black](7*0.8,0.7*-4)circle(0.1);
 			\filldraw[fill=black](11*0.8,0.7*0)circle(0.1);
 			\draw (8.5*0.8,0.7*0)--(7*0.8,0.7*-2);
 			\draw (8.5*0.8,0.7*0)--(10*0.8,0.7*-2);
 			\draw (7*0.8,0.7*-4)--(7*0.8,0.7*-2);
 			\coordinate [label={[text=black]above:$1$}] (x) at (8.5*0.8,0.7*0);
 			\coordinate [label={[text=black]above:$3$}] (x) at (11*0.8,0.7*0);
 			\coordinate [label={[text=black]below:$2$}] (x) at (7*0.8,0.7*-4);
 			
 			\filldraw[fill=black](15.5*0.8,0.7*0)circle(0.1);
 			\filldraw[fill=black](14*0.8,0.7*-2)circle(0.1);
 			\filldraw[fill=black](17*0.8,0.7*-2)circle(0.1);
 			\filldraw[fill=black](14*0.8,0.7*-4)circle(0.1);
 			\filldraw[fill=black](18*0.8,0.7*0)circle(0.1);
 			\draw (15.5*0.8,0.7*0)--(14*0.8,0.7*-2);
 			\draw (15.5*0.8,0.7*0)--(17*0.8,0.7*-2);
 			\draw (14*0.8,0.7*-4)--(14*0.8,0.7*-2);
 			\coordinate [label={[text=black]above:$1$}] (x) at (15.5*0.8,0.7*0);
 			\coordinate [label={[text=black]above:$2$}] (x) at (18*0.8,0.7*0);
 			\coordinate [label={[text=black]below:$3$}] (x) at (14*0.8,0.7*-4);
 			
 			\filldraw[fill=black](22.5*0.8,0.7*0)circle(0.1);
 			\filldraw[fill=black](21*0.8,0.7*-2)circle(0.1);
 			\filldraw[fill=black](24*0.8,0.7*-2)circle(0.1);
 			\filldraw[fill=black](24*0.8,0.7*-4)circle(0.1);
 			\filldraw[fill=black](25*0.8,0.7*0)circle(0.1);
 			\draw (22.5*0.8,0.7*0)--(21*0.8,0.7*-2);
 			\draw (22.5*0.8,0.7*0)--(24*0.8,0.7*-2);
 			\draw (24*0.8,0.7*-4)--(24*0.8,0.7*-2);
 			\coordinate [label={[text=black]above:$1$}] (x) at (22.5*0.8,0.7*0);
 			\coordinate [label={[text=black]above:$2$}] (x) at (25*0.8,0.7*0);
 			\coordinate [label={[text=black]below:$3$}] (x) at (24*0.8,0.7*-4);
 			
 			\coordinate [label=below:$F_1$] (x) at (2*0.8,0.7*-6);
 			\coordinate [label=below:$F_2$] (x) at (9*0.8,0.7*-6);
 			\coordinate [label=below:$F_3$] (x) at (16*0.8,0.7*-6);
 			\coordinate [label=below:$F_4$] (x) at (23*0.8,0.7*-6);
 			
 			\filldraw[fill=black](-1.5*0.8,0.7*-10)circle(0.1);
 			\filldraw[fill=black](-3*0.8,0.7*-12)circle(0.1);
 			\filldraw[fill=black](0*0.8,0.7*-12)circle(0.1);
 			\filldraw[fill=black](0*0.8,0.7*-14)circle(0.1);
 			\filldraw[fill=black](-1.5*0.8,0.7*-16)circle(0.1);
 			\filldraw[fill=black](-1.5*0.8,0.7*-18)circle(0.1);
 			\filldraw[fill=black](1.5*0.8,0.7*-16)circle(0.1);
 			\draw (-1.5*0.8,0.7*-10)--(0*0.8,0.7*-12);
 			\draw (-1.5*0.8,0.7*-10)--(-3*0.8,0.7*-12);
 			\draw (0*0.8,0.7*-14)--(0*0.8,0.7*-12);
 			\draw (-1.5*0.8,0.7*-16)--(0*0.8,0.7*-14);
 			\draw (1.5*0.8,0.7*-16)--(0*0.8,0.7*-14);
 			\draw (-1.5*0.8,0.7*-16)--(-1.5*0.8,0.7*-18);
 			\coordinate [label={[text=black]above:$1$}] (x) at(-1.5*0.8,0.7*-10);
 			\coordinate [label={[text=black]right:$2$}] (x) at(0*0.8,0.7*-14);
 			\coordinate [label={[text=black]below:$3$}] (x) at(-1.5*0.8,0.7*-18);
 			
 			\filldraw[fill=black](5*0.8,0.7*-10)circle(0.1);
 			\filldraw[fill=black](3.5*0.8,0.7*-12)circle(0.1);
 			\filldraw[fill=black](6.5*0.8,0.7*-12)circle(0.1);
 			\filldraw[fill=black](6.5*0.8,0.7*-14)circle(0.1);
 			\filldraw[fill=black](5*0.8,0.7*-16)circle(0.1);
 			\filldraw[fill=black](8*0.8,0.7*-18)circle(0.1);
 			\filldraw[fill=black](8*0.8,0.7*-16)circle(0.1);
 			\draw (5*0.8,0.7*-10)--(6.5*0.8,0.7*-12);
 			\draw (5*0.8,0.7*-10)--(3.5*0.8,0.7*-12);
 			\draw (6.5*0.8,0.7*-14)--(6.5*0.8,0.7*-12);
 			\draw (5*0.8,0.7*-16)--(6.5*0.8,0.7*-14);
 			\draw (8*0.8,0.7*-16)--(6.5*0.8,0.7*-14);
 			\draw (8*0.8,0.7*-16)--(8*0.8,0.7*-18);
 			\coordinate [label={[text=black]above:$1$}] (x) at(5*0.8,0.7*-10);
 			\coordinate [label={[text=black]right:$2$}] (x) at(6.5*0.8,0.7*-14);
 			\coordinate [label={[text=black]below:$3$}] (x) at(8*0.8,0.7*-18);
 			
 			\filldraw[fill=black](12*0.8,0.7*-10)circle(0.1);
 			\filldraw[fill=black](14.5*0.8,0.7*-10)circle(0.1);
 			\filldraw[fill=black](10.5*0.8,0.7*-12)circle(0.1);
 			\filldraw[fill=black](13.5*0.8,0.7*-12)circle(0.1);
 			\filldraw[fill=black](13.5*0.8,0.7*-14)circle(0.1);
 			\draw (12*0.8,0.7*-10)--(10.5*0.8,0.7*-12);
 			\draw (12*0.8,0.7*-10)--(13.5*0.8,0.7*-12);
 			\draw (13.5*0.8,0.7*-14)--(13.5*0.8,0.7*-12);
 			\coordinate [label={[text=black]above:$1$}] (x) at(12*0.8,0.7*-10);
 			\coordinate [label={[text=black]below:$2$}] (x) at(13.5*0.8,0.7*-14);
 			\coordinate [label={[text=black]above:$3$}] (x) at(14.5*0.8,0.7*-10);
 			
 			\filldraw[fill=black](19.5*0.8,0.7*-10)circle(0.1);
 			\filldraw[fill=black](17*0.8,0.7*-10)circle(0.1);
 			\filldraw[fill=black](18*0.8,0.7*-12)circle(0.1);
 			\filldraw[fill=black](21*0.8,0.7*-12)circle(0.1);
 			\filldraw[fill=black](21*0.8,0.7*-14)circle(0.1);
 			\draw (19.5*0.8,0.7*-10)--(18*0.8,0.7*-12);
 			\draw (19.5*0.8,0.7*-10)--(21*0.8,0.7*-12);
 			\draw (21*0.8,0.7*-14)--(21*0.8,0.7*-12);
 			\coordinate [label={[text=black]above:$2$}] (x) at(19.5*0.8,0.7*-10);
 			\coordinate [label={[text=black]below:$3$}] (x) at(21*0.8,0.7*-14);
 			\coordinate [label={[text=black]above:$1$}] (x) at(17*0.8,0.7*-10);
 			
 			\filldraw[fill=black](25*0.8,0.7*-10)circle(0.1);
 			\filldraw[fill=black](23.5*0.8,0.7*-12)circle(0.1);
 			\filldraw[fill=black](26.5*0.8,0.7*-12)circle(0.1);
 			\filldraw[fill=black](25*0.8,0.7*-14)circle(0.1);
 			\filldraw[fill=black](28*0.8,0.7*-14)circle(0.1);
 			\draw (25*0.8,0.7*-10)--(23.5*0.8,0.7*-12);
 			\draw (25*0.8,0.7*-10)--(26.5*0.8,0.7*-12);
 			\draw (26.5*0.8,0.7*-12)--(25*0.8,0.7*-14);
 			\draw (28*0.8,0.7*-14)--(26.5*0.8,0.7*-12);
 			\coordinate [label={[text=black]above:$1$}] (x) at(25*0.8,0.7*-10);
 			\coordinate [label={[text=black]below:$2$}] (x) at(25*0.8,0.7*-14);
 			\coordinate [label={[text=black]below:$3$}] (x) at(28.5*0.8,0.7*-14);
 			
 			\coordinate [label=below:$F_5$] (x) at (-1*0.8,0.7*-20);
 			\coordinate [label=below:$F_6$] (x) at (6*0.8,0.7*-20);
 			\coordinate [label=below:$F_7$] (x) at (13*0.8,0.7*-20);
 			\coordinate [label=below:$F_8$] (x) at (19*0.8,0.7*-20);
 			\coordinate [label=below:$F_9$] (x) at (26*0.8,0.7*-20);
 			
 		\end{tikzpicture}
 	\end{center}
 	\caption{All $9$ forests in  $\overline{\mathcal{F}}_{3}(2)$. }\label{overlineF}
 \end{figure}

\begin{definition}[Removable old leaf]\label{defioleaf}
	Given a forest $F=(T_1, T_2, \ldots, T_m)\in \mathcal{F}_{n}(k)$, let $u$ be an old   leaf of $T_i$ for some $i\in [m]$. Then the old leaf $u$ is said to be   removable   if  it verifies the following properties. 
\begin{itemize}
	\item The node $u$ is a grand child of the root of $T_i$, which is the only leaf  among all the grand children of the root of $T_i$;
	\item The first $k-1$ children of the root of $T_i$ are leaves;
 	\item The label of $u$ is smaller than that of the root of $T_{i+1}$ if $i<m$.
\end{itemize}
 For the forest  in Fig.~\ref{F}, the old leaf  $3$ is  removable, whereas the old leaves  $8$ and $10$ are not  removable. 
\end{definition}

\begin{figure}  
    \begin{center}
        \begin{tikzpicture}[font =\small , scale = 0.4]
        \filldraw[fill=black](0*0.8,0.7*0)circle(0.1);
        \filldraw[fill=black](3*0.8,0.7*0)circle(0.1);
        \filldraw[fill=black](1.5*0.8,0.7*-2)circle(0.1);
        \filldraw[fill=black](4.5*0.8,0.7*-2)circle(0.1);
        \filldraw[fill=black](1.5*0.8,0.7*-4)circle(0.1);
        \draw (3*0.8,0.7*0)--(1.5*0.8,0.7*-2);
        \draw (1.5*0.8,0.7*-4)--(1.5*0.8,0.7*-2);
        \draw (3*0.8,0.7*0)--(4.5*0.8,0.7*-2);
        \coordinate [label={[text=black]above:$1$}] (x) at(0*0.8,0.7*0);
        \coordinate [label={[text=black]above:$2$}] (x) at(3*0.8,0.7*0);
        \coordinate [label={[text=black]below:$3$}] (x) at(1.5*0.8,0.7*-4);

        \filldraw[fill=black](9.5*0.8,0.7*0)circle(0.1);
        \filldraw[fill=black](8*0.8,0.7*-2)circle(0.1);
        \filldraw[fill=black](11*0.8,0.7*-2)circle(0.1);
        \filldraw[fill=black](8*0.8,0.7*-4)circle(0.1);
        \filldraw[fill=black](6.5*0.8,0.7*-6)circle(0.1);
        \filldraw[fill=black](9.5*0.8,0.7*-6)circle(0.1);
        \filldraw[fill=black](6.5*0.8,0.7*-8)circle(0.1);
        \draw (9.5*0.8,0.7*0)--(8*0.8,0.7*-2);
        \draw (9.5*0.8,0.7*0)--(11*0.8,0.7*-2);
        \draw (8*0.8,0.7*-4)--(8*0.8,0.7*-2);
        \draw (9.5*0.8,0.7*-6)--(8*0.8,0.7*-4);
        \draw (6.5*0.8,0.7*-6)--(8*0.8,0.7*-4);
        \draw (6.5*0.8,0.7*-8)--(6.5*0.8,0.7*-6);
        \coordinate [label={[text=black]above:$1$}] (x) at(9.5*0.8,0.7*0);
        \coordinate [label={[text=black]right:$2$}] (x) at(8*0.8,0.7*-4);
        \coordinate [label={[text=black]below:$3$}] (x) at(6.5*0.8,0.7*-8);

        \filldraw[fill=black](16*0.8,0.7*0)circle(0.1);
        \filldraw[fill=black](14.5*0.8,0.7*-2)circle(0.1);
        \filldraw[fill=black](17.5*0.8,0.7*-2)circle(0.1);
        \filldraw[fill=black](14.5*0.8,0.7*-4)circle(0.1);
        \filldraw[fill=black](13*0.8,0.7*-6)circle(0.1);
        \filldraw[fill=black](16*0.8,0.7*-6)circle(0.1);
        \filldraw[fill=black](16*0.8,0.7*-8)circle(0.1);
        \draw (16*0.8,0.7*0)--(14.5*0.8,0.7*-2);
        \draw (16*0.8,0.7*0)--(17.5*0.8,0.7*-2);
        \draw (14.5*0.8,0.7*-4)--(14.5*0.8,0.7*-2);
        \draw (16*0.8,0.7*-6)--(14.5*0.8,0.7*-4);
        \draw (13*0.8,0.7*-6)--(14.5*0.8,0.7*-4);
        \draw (16*0.8,0.7*-8)--(16*0.8,0.7*-6);
        \coordinate [label={[text=black]above:$1$}] (x) at(16*0.8,0.7*0);
        \coordinate [label={[text=black]right:$2$}] (x) at(14.5*0.8,0.7*-4);
        \coordinate [label={[text=black]below:$3$}] (x) at(16*0.8,0.7*-8);

        \coordinate [label=below:$F_1$] (x) at (2.5*0.8,0.7*-10);
        \coordinate [label=below:$F_2$] (x) at (9*0.8,0.7*-10);
        \coordinate [label=below:$F_3$] (x) at (16*0.8,0.7*-10);

        \filldraw[fill=black](21.5*0.8,0.7*0)circle(0.1);
        \filldraw[fill=black](20*0.8,0.7*-2)circle(0.1);
        \filldraw[fill=black](23*0.8,0.7*-2)circle(0.1);
        \filldraw[fill=black](20*0.8,0.7*-4)circle(0.1);
        \filldraw[fill=black](23*0.8,0.7*-4)circle(0.1);
        \draw (21.5*0.8,0.7*0)--(20*0.8,0.7*-2);
        \draw (21.5*0.8,0.7*0)--(23*0.8,0.7*-2);
        \draw (20*0.8,0.7*-2)--(20*0.8,0.7*-4);
        \draw (23*0.8,0.7*-2)--(23*0.8,0.7*-4);
        \coordinate [label={[text=black]above:$1$}] (x) at(21.5*0.8,0.7*0);
        \coordinate [label={[text=black]below:$2$}] (x) at(23*0.8,0.7*-4);
        \coordinate [label={[text=black]below:$3$}] (x) at(20*0.8,0.7*-4);

        \filldraw[fill=black](35*0.8,0.7*0)circle(0.1);
        \filldraw[fill=black](33.5*0.8,0.7*-2)circle(0.1);
        \filldraw[fill=black](36.5*0.8,0.7*-2)circle(0.1);
        \filldraw[fill=black](33.5*0.8,0.7*-4)circle(0.1);
        \filldraw[fill=black](36.5*0.8,0.7*-4)circle(0.1);
        \draw (35*0.8,0.7*0)--(33.5*0.8,0.7*-2);
        \draw (35*0.8,0.7*0)--(36.5*0.8,0.7*-2);
        \draw (33.5*0.8,0.7*-2)--(33.5*0.8,0.7*-4);
        \draw (36.5*0.8,0.7*-2)--(36.5*0.8,0.7*-4);
        \coordinate [label={[text=black]above:$1$}] (x) at(35*0.8,0.7*0);
        \coordinate [label={[text=black]below:$3$}] (x) at(36.5*0.8,0.7*-4);
        \coordinate [label={[text=black]below:$2$}] (x) at(33.5*0.8,0.7*-4);

        \filldraw[fill=black](29*0.8,0.7*0)circle(0.1);
        \filldraw[fill=black](27.5*0.8,0.7*-2)circle(0.1);
        \filldraw[fill=black](30.5*0.8,0.7*-2)circle(0.1);
        \filldraw[fill=black](26*0.8,0.7*-4)circle(0.1);
        \filldraw[fill=black](29*0.8,0.7*-4)circle(0.1);
        \draw (29*0.8,0.7*0)--(30.5*0.8,0.7*-2);
        \draw (29*0.8,0.7*0)--(27.5*0.8,0.7*-2);
        \draw (27.5*0.8,0.7*-2)--(29*0.8,0.7*-4);
        \draw (27.5*0.8,0.7*-2)--(26*0.8,0.7*-4);
        \coordinate [label={[text=black]above:$1$}] (x) at(29*0.8,0.7*0);
        \coordinate [label={[text=black]below:$3$}] (x) at(29*0.8,0.7*-4);
        \coordinate [label={[text=black]below:$2$}] (x) at(26*0.8,0.7*-4);

        \coordinate [label=below:$F_4$] (x) at (21.5*0.8,0.7*-10);
        \coordinate [label=below:$F_5$] (x) at (28*0.8,0.7*-10);
        \coordinate [label=below:$F_6$] (x) at (35*0.8,0.7*-10);
         \end{tikzpicture}
    \end{center}
 		\caption{All $6$ forests of $\widehat{\mathcal{F}}_{3}(2)$. }\label{widehatF}
 \end{figure}
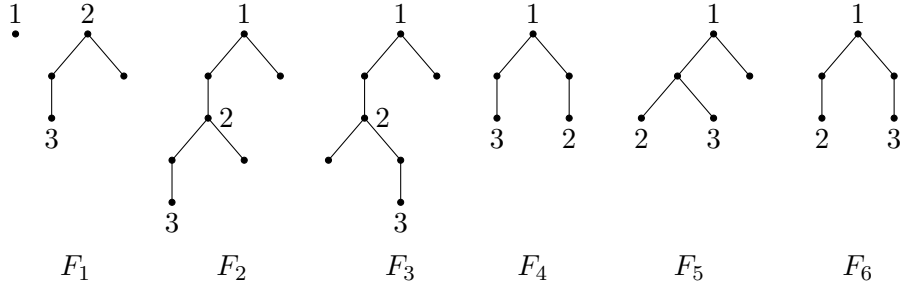

 Denote  by $\overline{\mathcal{F}}_{n}(k)$ the  subset of forests  $F\in \mathcal{F}_{n}(k)$ in which either the rightmost tree of $F$ is a singleton  or   the first $k-1$ children of the root in the rightmost tree of $F$ are leaves.   
 Let $\widehat{\mathcal{F}}_{n}(k)=\mathcal{F}_{n}(k)-\overline{\mathcal{F}}_{n}(k)$.
 For instance,  we draw all $9$ forests of $\overline{\mathcal{F}}_{3}(2)$  in Fig.~\ref{overlineF} and all $6$ forests of $\widehat{\mathcal{F}}_{3}(2)$ in Fig.~\ref{widehatF}.

 Now we are in the position to state our main result, which provides a  combinatorial interpretation for the   bi-$\gamma$-coefficients of $A^{(k)}_{n}(x)$.
  \begin{theorem}\label{gammaa}
 	  For $n, k\geq 1$, we have   
\begin{equation}\label{eqa}
 	  	a^{(k)}_{n}(x)=\sum\limits_{F\in \overline{\mathcal{F}}_{n}(k)}x^{\lleaf(F)-\si(F)}=\sum_{i=0}^{\lfloor{n-1\over 2}\rfloor} \overline{\gamma}_{n, k,i}x^i(1+x)^{n-1-2i}
 	  \end{equation}
	  and
	   \begin{equation}\label{eqb}
  xb^{(k)}_{n}(x)=\sum\limits_{F\in \widehat{\mathcal{F}}_{n}(k)}x^{\lleaf(F)-\si(F)}=\sum_{i=1}^{\lfloor{n\over 2}\rfloor} \widehat{\gamma}_{n,k, i}x^i(1+x)^{n-2i},
  	\end{equation}
 	where $\overline{\gamma}_{n, k,i}$ (resp.,~$\widehat{\gamma}_{n,k, i}$) enumerates the forests in $\overline{\mathcal{F}}_{n}(k)$ (resp.,~$\widehat{\mathcal{F}}_{n}(k)$) with $i$ old     leaves and  without any young   leaves or removable  old leaves. 
 \end{theorem}
 
 \begin{example} Take $n=3$ and $k=2$ in Theorem~\ref{gammaa}. The $\gamma$-positivity expansion of $a_3^{(2)}(x)$ is 
 $$
 a^{(2)}_3(x) = (1+x)^2+ 5x,
 $$
 where the $\gamma$-coefficients $1$ counts the forest $F_1$ and $5$ counts the forests $F_i$ for all $2\leq i\leq 6$ in Fig.~\ref{overlineF}. The $\gamma$-positivity expansion of $xb_3^{(2)}(x)$ is 
 $$
 xb^{(2)}_3(x) = 3x(1+x),
 $$
 where the $\gamma$-coefficient $3$ counts the forests $F_i$ for all $1\leq i\leq 3$ in Fig.~\ref{widehatF}. 
 \end{example}
  
The rest of this paper is devoted to the proof of Theorem~\ref{gammaa}.  Starting from the interpretations of $a^{(k)}_n(x)$ and $xb^{(k)}_n(x)$ in Proposition~\ref{prop2}, we construct in Section~\ref{Sec:3}  bijections between $k$-Stirling permutations and increasing pruned even $k$-ary forests, thereby proving the first equalities  stated in~\eqref{eqa} and~\eqref{eqb}. In Section~\ref{Sec:4}, we introduce a generalized Foata–Strehl  action on increasing pruned even $k$-ary trees which implies the  longest ascent-plateau polynomials  over  $k$-Stirling permutations with initial letter $1$ are $\gamma$-positive (see Theorem~\ref{gammac}). This $\gamma$-positivity result forms the second step of our proof of Theorem~\ref{gammaa} and may be of  independent interest. Finally in Section~\ref{Sec:5}, we develop two  transformations on increasing pruned even $k$-ary forests to conclude our proof of Theorem~\ref{gammaa}. This section constitutes the most technical  part of our combinatorial  approach.

\section{Bijections between permutations and forests}
\label{Sec:3}
 In this section, we   establish bijections between $k$-Stirling permutations and increasing pruned even $k$-ary forests, thereby proving the first equalities  stated in~\eqref{eqa} and~\eqref{eqb}.

 Let $k\geq 1$ and let $M=\{a_1, a_2, \ldots, a_n\}$ with $a_1<a_2<\ldots<a_n$.  Denote by $\mathcal{Q}_{M}(k)$ the set of $k$-Stirling permutations on $M_k$. 
  \begin{proposition}\label{propxi}
 	Let $M$ and $k$ be given as above.  
 There is a one-to-one correspondence $\xi$  between     $ \mathcal{Q}_{M}(k)$ and $\mathcal{F}_{M}(k)$ such that $\lap(\pi)=\lleaf(\xi(\pi))$ for any $\pi\in \mathcal{Q}_{M}(k)$.  
 \end{proposition}

\begin{figure}
     \begin{center}	\begin{tikzpicture}[font =\small , scale = 0.3]        \filldraw[fill=black](2,0)circle(0.1);        \filldraw[fill=black](2,-2)circle(0.1);        \filldraw[fill=black](0,-2)circle(0.1);        \filldraw[fill=black](4,-2)circle(0.1);        \filldraw[fill=black](0,-4)circle(0.1);        \filldraw[fill=black](4,-4)circle(0.1);        \draw(2,0)--(2,-2);        \draw(2,0)--(0,-2);        \draw(2,0)--(4,-2);        \draw(4,-4)--(2,-2);        \draw(0,-4)--(2,-2);        \coordinate [label=above:$1$] (x) at (2,0);        \coordinate [label=below:$3$] (x) at (0,-4);        \coordinate [label=below:$7$] (x) at (4,-4);        \filldraw[fill=black](9,0)circle(0.1);        \filldraw[fill=black](9,-2)circle(0.1);        \filldraw[fill=black](7,-2)circle(0.1);        \filldraw[fill=black](11,-2)circle(0.1);        \filldraw[fill=black](7,-4)circle(0.1);        \filldraw[fill=black](11,-4)circle(0.1);        \filldraw[fill=black](5,-6)circle(0.1);        \filldraw[fill=black](7,-6)circle(0.1);        \filldraw[fill=black](9,-6)circle(0.1);        \filldraw[fill=black](9,-8)circle(0.1);        \draw(9,0)--(9,-2);        \draw(9,0)--(7,-2);        \draw(9,0)--(11,-2);        \draw(11,-4)--(11,-2);        \draw(7,-4)--(7,-2);        \draw(5,-6)--(7,-4);        \draw(9,-6)--(7,-4);        \draw(7,-6)--(7,-4);        \draw(9,-6)--(9,-8);        \coordinate [label=above:$2$] (x) at (9,0);        \coordinate [label=left:$4$] (x) at (7,-4);        \coordinate [label=below:$5$] (x) at (11,-4);        \coordinate [label=below:$6$] (x) at (9,-8);        \filldraw[fill=black](14,0)circle(0.1);        \coordinate [label=above:$8$] (x) at (14,0);        \coordinate [label=below:$ T_1$] (x) at (2,-10);        \coordinate [label=below:$ T_2$] (x) at (9,-10);        \coordinate [label=below:$ T_3$] (x) at (14,-10);                \end{tikzpicture}    \end{center}
  \caption{An increasing pruned even $3$-ary forest $(T_1, T_2, T_3)$. }\label{xi}
\end{figure}
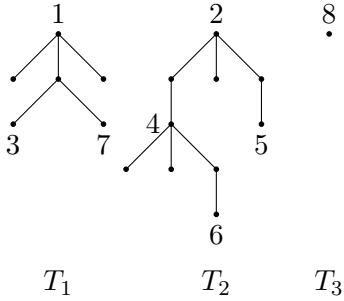

\pf Given a $k$-Stirling permutation $\pi\in \mathcal{Q}_{M}(k)$,  if $\pi$ is empty, set $\xi(\pi)$ to be the empty forest. 
Otherwise, suppose that $\pi$ has $m$ right-to-left minima, say $b_1, b_2, \ldots, b_m$ with $b_1<b_2<\cdots<b_m$. 
Recall that    a  right-to-left minimum of the  permutation  $\sigma=\sigma_1\sigma_2\ldots \sigma_n$ is a letter $\sigma_i$ such that  $\sigma_j>\sigma_i$ for all $j>i$.  Then $\pi$ can be uniquely decomposed as $w_1w_{2}\ldots w_m$, where each $w_i$ is a $k$-Stirling permutation that ends with $b_i$. Clearly, each $w_i$ ends with the smallest element of $w_i$.  For each $i\in [m]$, generate a tree $T_i$ as follows:
\begin{itemize}
	\item Decompose $w_i$ uniquely as
	$\mu_1b_i\mu_2 b_i\ldots \mu_kb_i$ where each $\mu_i$ is either empty or a  $k$-Stirling permutation.
	\item 	 If   $\mu_j$ is empty for all $j\in [k]$, set $T_i$ to be  a singleton with label $b_i$;
		\item Otherwise, set the node with label $b_i$ to be the root of $T_i$. Then attach $k$ unlabeled leaves $v_1, v_2, \ldots, v_k$ from left to right to the root,  and attach all the trees in the forest $\xi(\mu_j)$ to the node $v_j$ as subtrees from left to right for all $j\in [k]$. Denote by $T_i$ the resulting tree.
	\end{itemize}
Set $\xi(\pi)=(T_1, T_2, \ldots, T_m)$. For example, if we  let $\pi=133377711446664225552888\in \mathcal{Q}_{8}(3)$,  then its corresponding forest $\xi(\pi)=(T_1, T_2, T_3)$ is illustrated in Fig.~\ref{xi}. 

Now we proceed to show that $\xi(\pi)\in \mathcal{F}_{M}(k)$ and  $\lap(\pi)=\lleaf(\xi(\pi))$ by induction on the $|M|$. 
Clearly, when $|M|=1$, the resulting forest $\xi(\pi)$ is a singleton. In this case, we have $\xi(\pi)\in \mathcal{F}_{M}(k)$ and $\lap(\pi)=\lleaf(\xi(\pi))=1$. Now we assume that $\xi(\pi)\in \mathcal{F}_{N}(k)$ and  $\lap(\pi)=\lleaf(\xi(\pi))$ for any permutation $\pi\in \mathcal{Q}_{N}(k)$   with $|N|<|M|$. Observe that $\lap(\pi)=\sum\limits_{i=1}^{m}\lap(w_i)$ and $\lleaf(T)=\sum\limits_{i=1}^{m}\lleaf(T_i)$.
  In order to show that  $\xi(\pi)\in \mathcal{F}_{M}(k)$ and $\lap(\pi)=\lleaf(\xi(\pi))$, it suffices to show that each $T_i$ is an increasing pruned even $k$-ary tree  with $\lap(w_i)=\lleaf(T_i)$ for all $i\in [m]$. This can be justified by the construction of $T_i$ and the induction hypothesis. 
It is apparent that the construction of $\xi$ is reversible and hence  the map 
$\xi$ is the desired   bijection. This completes the proof. \qed

Let $\widetilde{\mathcal{Q}}_{M}(k)$  be the subset of $k$-Stirling permutations in $\mathcal{Q}_{M}(k)$ which starts with the smallest element of $M$.  
\begin{proposition}\label{propchi}
Let $M$ and $k$ be given as above. For  $n\geq 2$, there is a one-to-one correspondence $\chi$ between $\widetilde{\mathcal{Q}}_{M}(k)$ and $\mathcal{T}_{M}(k)$ such that $\ap(\pi)=\lleaf(\chi(\pi))$ for any $\pi\in \widetilde{\mathcal{Q}}_{M}(k)$. Moreover, for any $\pi\in \widetilde{\mathcal{Q}}_{M}(k)$, the permutation $\pi$ starts with a longest  plateau if and only if  the first $k-1$ children of the root of the corresponding tree are leaves. 
\end{proposition}
\pf  Recall that $a_1$ is the smallest element of $M$.  Given a $k$-Stirling permutation $\pi\in \widetilde{\mathcal{Q}}_{M}(k)$,   then $\pi$ can be uniquely decomposed as $ a_1 \, w_1 \, a_1 \, \cdots  \, w_{k-1} \,
a_1 \, w_k$, where $w_j$ is either empty or a $k$-Stirling permutation for all $1 \leq j\leq k$.
Set the node with label  $a_1$ to be the root of $\chi(\pi)$ and   attach $k$ unlabeled leaves to the root, say $v_1, v_2, \ldots, v_k$ from left to right. For all  $j\in[k]$, if $w_j$ is non-empty, attach all the trees in the forest  $\xi(w_j)$ to the node $v_j$ as subtrees  and arrange the subtrees of $v_j$ so that the children of the node $v_j$ are increasing from left to right. For instance, if we let $\pi=244666422555$, then its corresponding tree  $\chi(\pi)=T_2$ is illustrated in Fig.~\ref{xi}.

 By Proposition $\ref{propxi}$, the resulting tree $\chi(\pi)$ is a tree in  $\mathcal{T}_{M}(k)$ such that $\ap(\pi)=\lleaf(\chi(\pi))$.  Moreover,  the permutation $\pi$ starts with a  longest plateau if and only if  the first $k-1$ children of the root of  $\chi(\pi)$ are leaves. 
It is apparent that the construction of $\chi$ is reversible and hence  the map 
$\chi$ is the desired   bijection. This completes the proof. \qed

\begin{proposition}\label{propzeta}
	For $n,k\geq 1$, there is a one-to-one correspondence $\zeta$ between $ \mathcal{Q}_{n}(k)$ and $\mathcal{F}_{n}(k)$ such that $\ap(\pi)=\lleaf(\zeta(\pi))-\si(\zeta(\pi))$ for any $\pi\in \mathcal{Q}_{n}(k)$. Moreover, we have $\pi\in \overline{\mathcal{Q}}_{n}(k) $ if and only if  $\zeta(\pi)\in \overline{\mathcal{F}}_{n}(k)$.
\end{proposition}
\pf Given a $k$-Stirling permutation $\pi \in  \mathcal{Q}_{n}(k)$,  assume that $\pi$ has $m$ left-to-right minima, say   $a_1, a_2, \ldots, a_m$ with $a_1<a_2<\ldots<a_m$. Recall that    a  left-to-right minimum of the  permutation  $\sigma=\sigma_1\sigma_2\ldots \sigma_n$ is a letter $\sigma_i$ such that  $\sigma_j>\sigma_i$ for all $j<i$.  Then $\pi$ can be uniquely decomposed as $w_mw_{m-1}\ldots w_1$, where each $w_i$ is a $k$-Stirling permutation that starts with $a_i$. Clearly, each $w_i$ starts with the smallest element of $w_i$. Then set $\zeta(\pi)$ to be the forest $(T_1, T_2, \ldots, T_m)$ where 
\begin{itemize}
	\item $T_i$ has only one node in which  the root is labeled by $a_i$ if $w_i$ has only one distinct element;
	\item otherwise, $T_i=\chi(w_i)$. 
\end{itemize}
For example,  if we let $\pi=888244666422555113337771\in \mathcal{Q}_8(3)$, then its corresponding forest $\zeta(\pi)=(T_1, T_2, T_3)$ is illustrated in Fig.~\ref{xi}.  
By Proposition \ref{propchi}, one can easily check that the resulting forest $\zeta(\pi)\in \mathcal{F}_n(k)$ and $\pi\in \overline{\mathcal{Q}}_{n}(k) $ if and only if  $\zeta(\pi)\in \overline{\mathcal{F}}_{n}(k)$. 
Again by Proposition \ref{propchi}, we have 
$\ap(w_i)=\lleaf(T_i)$ if $w_i$ has at least two distinct elements.  By the construction of $T_i$,  the tree $T_i$ is singleton and $\lleaf(T_i)=1$ if  $w_i$ has only one  distinct element.  However, we have $\ap(w_i)=0$ when  $w_i$ has only one  distinct element.   Hence, we deduce that $\ap(\pi)=\sum\limits_{i=1}^{m}\ap(w_i)=\lleaf(\zeta(\pi))-\si(\zeta(\pi))$ as desired.
Since $\chi$ is a bijection by Proposition \ref{propchi}, the construction of the map $\zeta$ is reversible and hence the map $\zeta$ is a bijection, completing the proof. 
 \qed
 
It then follows from Propositions~\ref{propzeta} and~\ref{prop2}  that
 \begin{equation}\label{eq:ab}
 a^{(k)}_{n}(x)=\sum\limits_{F\in \overline{\mathcal{F}}_{n}(k)}x^{\lleaf(F)-\si(F)}\quad\text{and}\quad xb^{(k)}_{n}(x)=\sum\limits_{F\in \widehat{\mathcal{F}}_{n}(k)}x^{\lleaf(F)-\si(F)},
 \end{equation}
 proving the first equalities  stated in~\eqref{eqa} and~\eqref{eqb}.
 
 \section{A generalized Foata–Strehl   action on  increasing pruned even $k$-ary trees } 
 \label{Sec:4}
 In this section, we develop a generalized Foata–Strehl   action on  increasing pruned even $k$-ary trees     in the sprit of  the generalized Foata-Strehl
action  on weakly increasing trees \cite{Lin}.

For $n,k\geq 1$, 
define 
$$
\widetilde{\mathcal{Q}}_{n}(k)=\{\pi=\pi_1\pi_2\ldots \pi_{kn}\in \mathcal{Q}_n(k)\mid \pi_1=1 \},
$$
 that is $\widetilde{\mathcal{Q}}_{n}(k)=\widetilde{\mathcal{Q}}_{[n]}(k)$. 
Define the polynomial  $c^{(k)}_{n}(x)$ by $$ c^{(k)}_{n}(x)=\sum\limits_{\pi\in \widetilde{\mathcal{Q}}_{n}(k)}x^{\ap(\pi)}.$$ 
It is clear that
 \begin{equation}\label{reca}
c_{n+1}^{(1)}(x)=xA_n(x).
 \end{equation}
 In \cite{Lin}, Lin-Ma-Ma-Zhou provided  an increasing  tree  interpretation of the $\gamma$-coefficients of the Eulerian polynomials $A_n(x)$.
In analogy to   the Eulerian polynomials, we obtain the following $\gamma$-expansion  of the  polynomials $c_n^{(k)}(x)$. 
 \begin{theorem}\label{gammac}
Set $n\geq 2$ and  $k\geq 1$. Then we have 
\begin{equation}\label{eqc}
 c^{(k)}_{n}(x)=\sum\limits_{T\in \mathcal{T}_{n}(k)}x^{\lleaf(T)}=\sum_{i=1}^{\lfloor{n\over 2}\rfloor} \widetilde{\gamma}_{n,k, i}x^i(1+x)^{n-2i},
 		\end{equation}
 	where
 	$ \widetilde{\gamma}_{n,k, i}$ enumerates the trees in $\mathcal{T}_{n}(k)$ with $i$ labeled     leaves and without any young  leaves. 
 \end{theorem}
 
 The first few $\gamma$-positivity expansions of $	c^{(3)}_{n}(x)$ read as follows:
 $$
c^{(3)}_2(x) = 3x, \quad 	c^{(3)}_3(x)= 9x(1+x)\quad\text{and}\quad c^{(3)}_4(x)=27x(1+x)^2+54x^2.
 $$
  By contracting all the edges connecting the nodes on even levels and the nodes  on odd levels for any $T\in\mathcal{T}_{n}(1)$, we get an increasing tree $T'$. Hence,  in view of~\eqref{reca}, setting $k=1$ in~\eqref{eqc} recovers  the $\gamma$-positivity expansion of the Eulerian polynomial $A_n(x)$ proved in~\cite{Lin}.

Let $M=\{a_1, a_2, \ldots, a_n\}$ and $k\geq 1$. For $T\in \mathcal{T}_{M}(k)$, suppose that the node $v$ is labeled by $x$ for some $x\in M$.
 We define the  increasing pruned even $k$-ary tree $\Phi_x(T)$ as follows:
\begin{itemize}
\item  If $v$ is an old   internal node, then suppose that the node $u$ is the grand parent of $v$,  the nodes $\mu_1, \mu_2, \ldots, \mu_k$ are all the children of $u$ from left to right, and the nodes $\nu_1, \nu_2, \ldots, \nu_k$ are all the children of $v$ from left to right. For each $j\in [k]$,  delete the edge connecting the node $\nu_j$ and the node $v$, attach the subtrees of   $\nu_j$ to  
$\mu_j$ as subtrees  in the order that  the children of $\mu_j$ are increasing from left to right.  Denote by $\Phi_x(T)$ the resulting tree. 
\item If $v$ is a young   leave, then suppose that  node $u$ is the grand parent of $v$,  the nodes $\mu_1, \mu_2, \ldots, \mu_k$ are all the children of $u$ from left to right.  For all $j\in [k]$, let $F_j$  be the set of all the subtrees of the node $\mu_j$ whose roots  have labels greater than the label of the node $v$. First, we attach $k$ unlabeled leaves to the node of $v$, say $\nu_1, \nu_2, \ldots, \nu_k$. Then remove all the subtrees in $F_j$ from  $\mu_j$ and attach them to the the node $\nu_j$ as subtrees in the order that  the children of $\nu_j$ are increasing from left to right for all $j\in [k]$.  Denote by $\Phi_x(T)$ the resulting tree. 
\item If $v$ is neither an old   internal node nor a young   leaf, then let $\Phi_x(T)=T$. 
\end{itemize}
 
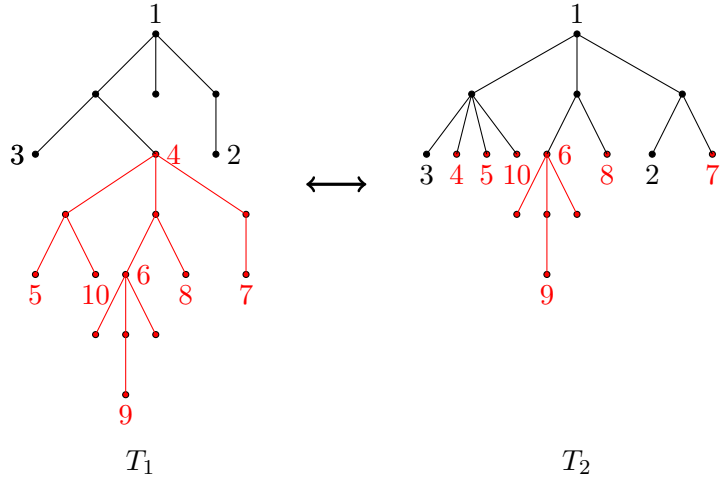
\begin{figure}    	\begin{center}	\begin{tikzpicture}[font =\small , scale = 0.4][<->,>=stealth,thick]		\filldraw[fill=black](4,0)circle(0.1);			\coordinate [label=above:$1$] (x) at (4,0);            \filldraw[fill=black](2,-2)circle(0.1);            \filldraw[fill=black](4,-2)circle(0.1);            \filldraw[fill=black](6,-2)circle(0.1);            \filldraw[fill=black](0,-4)circle(0.1);			\coordinate [label=left:$3$] (x) at (0,-4);            \filldraw[fill=red](4,-4)circle(0.1);            \filldraw[fill=black](6,-4)circle(0.1);            \coordinate [label={[text=red]right:$4$}] (x) at (4,-4);            \coordinate [label=left:$3$] (x) at (0,-4);            \coordinate [label=right:$2$] (x) at (6,-4);            \draw(2,-2)--(4,0);            \draw(4,-2)--(4,0);            \draw(6,-2)--(4,0);            \draw(0,-4)--(2,-2);            \draw(4,-4)--(2,-2);            \draw(6,-4)--(6,-2);            \filldraw[fill=red](1,-6)circle(0.1);            \filldraw[fill=red](4,-6)circle(0.1);            \filldraw[fill=red](7,-6)circle(0.1);            \filldraw[fill=red](0,-8)circle(0.1);            \filldraw[fill=red](7,-8)circle(0.1);            \filldraw[fill=red](5,-8)circle(0.1);            \filldraw[fill=red](3,-8)circle(0.1);            \filldraw[fill=red](2,-8)circle(0.1);            \draw[red] (1,-6)--(4,-4);             \draw[red] (4,-6)--(4,-4);            \draw[red] (7,-6)--(4,-4);            \draw[red] (1,-6)--(0,-8);            \draw[red] (1,-6)--(2,-8);            \draw[red] (4,-6)--(3,-8);            \draw[red] (4,-6)--(5,-8);            \draw[red] (7,-8)--(7,-6);            \coordinate [label={[text=red]below:$5$}] (x) at (0,-8);            \coordinate [label={[text=red]below:$10$}] (x) at (2,-8);            \coordinate [label={[text=red]right:$6$}] (x) at (3,-8);            \coordinate [label={[text=red]below:$8$}] (x) at (5,-8);            \coordinate [label={[text=red]below:$7$}] (x) at (7,-8);            \filldraw[fill=red](2,-10)circle(0.1);            \filldraw[fill=red](3,-10)circle(0.1);            \filldraw[fill=red](4,-10)circle(0.1);            \filldraw[fill=red](3,-12)circle(0.1);            \draw[red] (3,-8)--(4,-10);             \draw[red] (3,-8)--(3,-10);            \draw[red] (3,-8)--(2,-10);            \draw[red] (3,-12)--(3,-10);            \coordinate [label={[text=red]below:$9$}] (x) at (3,-12);            \draw[line width=1.2pt][<->] (9,-5) -- (11,-5);              \filldraw[fill=black](18,0)circle(0.1);             \filldraw[fill=black](14.5,-2)circle(0.1);             \filldraw[fill=black](21.5,-2)circle(0.1);             \filldraw[fill=black](13,-4)circle(0.1);             \filldraw[fill=red](14,-4)circle(0.1);             \filldraw[fill=red](15,-4)circle(0.1);             \filldraw[fill=red](16,-4)circle(0.1);             \filldraw[fill=red](17,-4)circle(0.1);             \filldraw[fill=red](19,-4)circle(0.1);             \filldraw[fill=black](20.5,-4)circle(0.1);             \filldraw[fill=red](22.5,-4)circle(0.1);             \filldraw[fill=black](18,-2)circle(0.1);             \draw(14.5,-2)--(18,0);             \draw(18,-2)--(18,0);             \draw(21.5,-2)--(18,0);             \draw(13,-4)--(14.5,-2);             \draw(14,-4)--(14.5,-2);             \draw(15,-4)--(14.5,-2);             \draw(16,-4)--(14.5,-2);             \draw(17,-4)--(18,-2);             \draw(19,-4)--(18,-2);             \draw(20.5,-4)--(21.5,-2);             \draw(22.5,-4)--(21.5,-2);             \coordinate [label=above:$1$] (x) at (18,0);             \coordinate [label=below:$3$] (x) at (13,-4);             \coordinate [label={[text=red]below:$10$}] (x) at (16,-4);             \coordinate [label={[text=red]below:$4$}] (x) at (14,-4);             \coordinate [label={[text=red]below:$5$}] (x) at (15,-4);             \coordinate [label={[text=red]right:$6$}] (x) at (17,-4);             \coordinate [label={[text=red]below:$8$}] (x) at (19,-4);             \coordinate [label={[text=black]below:$2$}] (x) at (20.5,-4);             \coordinate [label={[text=red]below:$7$}] (x) at (22.5,-4);             \filldraw[fill=red](16,-6)circle(0.1);             \filldraw[fill=red](17,-6)circle(0.1);             \filldraw[fill=red](18,-6)circle(0.1);             \filldraw[fill=red](17,-8)circle(0.1);             \draw[red](18,-6)--(17,-4);             \draw[red](17,-6)--(17,-4);             \draw[red](16,-6)--(17,-4);             \draw[red](17,-8)--(17,-6);             \coordinate [label={[text=red]below:$9$}] (x) at (17,-8);              \coordinate [label=above:$  T_1$] (x) at (3.5,-15);               \coordinate [label=above:$ T_2$] (x) at (18,-15);                \end{tikzpicture}             
 \end{center}
 \caption{ The  transformation $\Phi_{4}$. }\label{GFS1}
\end{figure}

 For example, let $T_1$ and $T_2$ be  increasing pruned even $3$-ary trees as shown in Fig.~\ref{GFS1}.  Then we have
 $\Phi_4(T_1)=T_2$ and $\Phi_4(T_2)=T_1$. 
 The following lemma  is clear from the definition of $\Phi_x(T)$.
 \begin{lemma}\label{lemPhi}
 	Let $T\in \mathcal{T}_{M}(k)$.  Then for any $x\in M$, we have $\Phi_x(T)\in \mathcal{T}_{M}(k)$. Moreover, the resulting tree $\Phi_x(T)$ verifies the following properties.
 	\begin{itemize} 
 		 	  		\item The node  with label  $x$ is an old  internal node of $T$ if and only if the node with label $x$ is a young   leaf of $\Phi_x(T)$;
 		\item The transformation  $\Phi_x$ preserves the types (old/young leaf/ internal   node ) of the other labeled  nodes. 
 	\end{itemize}	
 	
 \end{lemma}
In view of Lemma \ref{lemPhi}, it is not hard to see that the transformation $\Phi_x$ is an involution on  $\mathcal{T}_{M}(k)$ and $\Phi_x$ and $\Phi_y$ commute for all $x,y\in M$.  For any  $S\subseteq M$, we can define the function $\Phi_{S}: \mathcal{T}_{\mathcal{M}}(k) \rightarrow \mathcal{T}_{\mathcal{M}}(k) $ by $\Phi_{S}=\prod\limits_{x\in S}\Phi_x$, where the product is the functional compositions.  Hence the group  $\mathbb{Z}_2^{M}$ acts on  $\mathcal{T}_{\mathcal{M}}(k)$ via the function $\Phi_S$ for $S\subseteq M$. Such   an action is called the 
{\em generalized Foata–Strehl action} (GFS-action for short) on  increasing even $k$-ary trees.

Given a tree $T\in \mathcal{T}_{M}(k)$, let $\oleaf(T)$  (resp. $\oint(T)$, $\yleaf(T)$ ) denote the number of old  leaves (resp. old   internal nodes,  young   leaves)of $T$. We shall record the   the relation between $\oint(T)$ and $\oleaf(T)$ for any $T$ without any young leaves in the following lemma, which will play important roles in the proof of Theorem \ref{gammac}. 
\begin{lemma}\label{relation1}
Let $M$ and $k$ be given as above and let $n\geq 2$. For any tree $T\in \mathcal{T}_{M}(k)$ without any young leaves,   we have \begin{equation}\label{re3}
	\oint(T)=n-2\oleaf(T).\end{equation}		
\end{lemma}
\pf  
Clearly, for any tree $T\in \mathcal{T}_{M}(k) $, it is easily seen that
\begin{equation}\label{eq1}
	\oleaf(T)+\yleaf(T)+\lint(T)=n,
\end{equation}
where $\lint(T)$ denotes the number of labeled internal nodes of $T$. 
Since each labeled internal node has either  an old   leaf or an old   internal node as its grand child, we derive that
\begin{equation}\label{eq2}
	\oleaf(T)+\oint(T)=\lint(T).
\end{equation}
In view of (\ref{eq1}) and (\ref{eq2}),  we deduce that $$
	\oint(T)=n-2\oleaf(T)$$ since $\yleaf(T)=0$, completing the proof. \qed

We are ready for the proof of Theorem  \ref{gammac}.

  \noindent{\bf Proof of Theorem  \ref{gammac}.}
We aim  to prove the following $\gamma$-expansion 
\begin{equation}\label{gammaT}
\sum\limits_{T\in \mathcal{T}_M(k)}x^{\lleaf(T)}= \sum_{i=1}^{\lfloor{n\over 2}\rfloor}  \widetilde\gamma_{M, k, i}x^i(1+x)^{n-2i},
\end{equation}
where
$\widetilde\gamma_{M, k, i}$ enumerates the trees in $\mathcal{T}_{M}(k)$ with $i$   labeled  leaves and without any young   leaves. In view of  Proposition~\ref{propchi}, setting $M=[n]$ in~\eqref{gammaT}  leds to the proof of Theorem~\ref{gammac}.

  For each $T\in \mathcal{T}_{M}(k)$, we define $[T]$ to be the
set of  trees $S$ that can be transformed  to $T$
via the GFS-action on increasing pruned even $k$-ary trees. Clearly, the GFS-action divides $\mathcal{T}_{M}(k)$ into disjoint orbits. 
In view of Lemma  \ref{lemPhi}, each $[T]$ contains exactly one tree without any young   leaves that we denote by $\widehat{T}$.
Then, by Lemma~\ref{lemPhi}, we have
$$
\sum_{T\in [\widehat{T}] 
 }x^{\lleaf(T)} 
=  x^{\oleaf(\widehat{T})}(1+x)^{\oint(\widehat{T})}=x^{\oleaf(\widehat{T})}(1+x)^{n-2\oleaf(\widehat{T})},
$$
where the second equality follows from (\ref{re3}). 
Let $\mathcal{T}^*_{M}(k)$ denote the subset of trees without any young   leaves in $\mathcal{T}_M(k)$.   Then,  we derive that
$$
\begin{array}{lll}
	\sum\limits_{T\in \mathcal{T}_M(k)}x^{\lleaf(T)}&=&  \sum\limits_{\widehat{T}\in \mathcal{T}^*_{M}(k) }\sum\limits_{T\in [\widehat{T}] }x^{\lleaf(T)} \\
	&=&  \sum\limits_{\widehat{T}\in \mathcal{T}^*_{M}(k) }x^{\oleaf(\widehat{T})}(1+x)^{n-2\oleaf(\widehat{T})}, \\
	 		&=&  \sum\limits_{i=1}^{\lfloor{n\over 2}\rfloor}  \widetilde\gamma_{M,k,i} x^{i}(1+x)^{n-2i}, \\
\end{array}
$$ 
 completing the proof. \qed

\section{Proof of Theorem~\ref{gammaa}} 

\label{Sec:5}
 This section is devoted to finishing  the proof of Theorem~\ref{gammaa}. To this end, we  develop two transformations on increasing pruned even $k$-ary forests. 
 
 \subsection{GFS-action on $\mathcal{F}_n(k)$ and some basic properties}

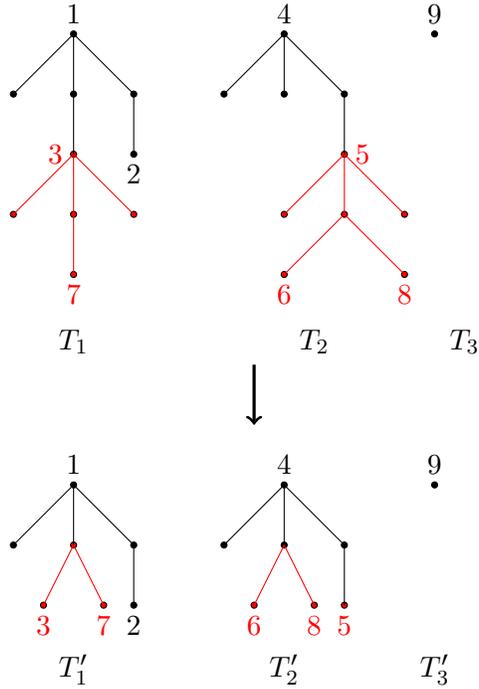
\begin{figure}    \begin{center}	\begin{tikzpicture}[font =\small , scale = 0.4]        \filldraw[fill=black](2,0)circle(0.1);        \filldraw[fill=black](0,-2)circle(0.1);        \filldraw[fill=black](2,-2)circle(0.1);        \filldraw[fill=black](4,-2)circle(0.1);        \filldraw[fill=black](4,-4)circle(0.1);        \filldraw[fill=red](2,-4)circle(0.1);        \filldraw[fill=red](0,-6)circle(0.1);        \filldraw[fill=red](2,-6)circle(0.1);        \filldraw[fill=red](4,-6)circle(0.1);        \filldraw[fill=red](2,-8)circle(0.1);        \draw[red] (2,-4)--(2,-6);        \draw[red] (2,-4)--(0,-6);        \draw[red] (2,-4)--(4,-6);        \draw[red] (2,-6)--(2,-8);        \draw (2,0)--(2,-2);        \draw (2,0)--(0,-2);        \draw (2,0)--(4,-2);        \draw (2,-4)--(2,-2);         \draw (4,-2)--(4,-4);        \coordinate [label={[text=red]below:$7$}] (x) at (2,-8);        \coordinate [label={[text=red]left:$3$}] (x) at (2,-4);        \coordinate [label={[text=black]below:$2$}] (x) at (4,-4);               \filldraw[fill=black](9,0)circle(0.1);        \filldraw[fill=black](14,0)circle(0.1);        \filldraw[fill=black](7,-2)circle(0.1);        \filldraw[fill=black](9,-2)circle(0.1);        \filldraw[fill=black](11,-2)circle(0.1);        \filldraw[fill=red](11,-4)circle(0.1);        \filldraw[fill=red](9,-6)circle(0.1);        \filldraw[fill=red](11,-6)circle(0.1);        \filldraw[fill=red](13,-6)circle(0.1);        \filldraw[fill=red](13,-8)circle(0.1);        \filldraw[fill=red](9,-8)circle(0.1);        \draw (7,-2)--(9,0);        \draw (9,-2)--(9,0);        \draw (11,-2)--(9,0);        \draw (11,-2)--(11,-4);        \draw[red] (9,-6)--(11,-4);        \draw[red] (11,-6)--(11,-4);        \draw[red] (13,-6)--(11,-4);        \draw[red] (11,-6)--(9,-8);        \draw[red] (11,-6)--(13,-8);        \coordinate [label={[text=black]above:$4$}] (x) at (9,0);        \coordinate [label={[text=black]above:$1$}] (x) at (2,0);        \coordinate [label={[text=red]right:$5$}] (x) at (11,-4);        \coordinate [label={[text=red]below:$6$}] (x) at (9,-8);        \coordinate [label={[text=red]below:$8$}] (x) at (13,-8);        \coordinate [label={[text=black]above:$9$}] (x) at (14,0);         \draw[line width=1.2pt][->] (8,-11) -- (8,-13);        \filldraw[fill=black](2,-15)circle(0.1);        \filldraw[fill=black](0,-17)circle(0.1);        \filldraw[fill=red](2,-17)circle(0.1);        \filldraw[fill=black](4,-17)circle(0.1);        \filldraw[fill=black](4,-19)circle(0.1);        \filldraw[fill=red](1,-19)circle(0.1);        \filldraw[fill=red](3,-19)circle(0.1);        \draw (2,-15)--(2,-17);        \draw (2,-15)--(0,-17);        \draw (2,-15)--(4,-17);        \draw (4,-17)--(4,-19);        \draw[red] (2,-17)--(1,-19);        \draw[red] (2,-17)--(3,-19);        \coordinate [label={[text=black]above:$1$}] (x) at (2,-15);        \coordinate [label={[text=red]below:$3$}] (x) at (1,-19);        \coordinate [label={[text=red]below:$7$}] (x) at (3,-19);         \coordinate [label={[text=black]below:$2$}] (x) at (4,-19);         \filldraw[fill=black](9,-15)circle(0.1);        \filldraw[fill=black](9,-17)circle(0.1);        \filldraw[fill=black](7,-17)circle(0.1);        \filldraw[fill=black](11,-17)circle(0.1);        \filldraw[fill=red](8,-19)circle(0.1);        \filldraw[fill=red](11,-19)circle(0.1);        \filldraw[fill=red](10,-19)circle(0.1);        \filldraw[fill=black](14,-15)circle(0.1);        \draw (9,-15)--(9,-17);        \draw (9,-15)--(7,-17);        \draw (9,-15)--(11,-17);        \draw (11,-19)--(11,-17);        \draw[red] (9,-17)--(8,-19);        \draw[red] (9,-17)--(10,-19);        \coordinate [label={[text=black]above:$9$}] (x) at (14,-15);        \coordinate [label={[text=black]above:$4$}] (x) at (9,-15);         \coordinate [label={[text=red]below:$6$}] (x) at (8,-19);         \coordinate [label={[text=red]below:$8$}] (x) at (10,-19);         \coordinate [label={[text=red]below:$5$}] (x) at (11,-19);          \coordinate [label=above:$ T_1$] (x) at (2,-11);        \coordinate [label=above:$ T_2$] (x) at (10,-11);        \coordinate [label=above:$ T_3$] (x) at (15,-11);        \coordinate [label=above:$ T'_1$] (x) at (2,-22);        \coordinate [label=above:$ T'_2$] (x) at (9,-22);        \coordinate [label=above:$ T'_3$] (x) at (14,-22);                        \end{tikzpicture}    \end{center}
		\caption{An example of the  transformation  $\Phi_{S}$ with $S=\{3,5\}$. }\label{PhiF}
\end{figure}

Given a forest $F=(T_1, T_2, \ldots, T_m)\in\mathcal{F}_n(k)$,  let  $\Oint(F)$ (resp.,  $\Si(F)$, $\Yleaf(F)$) denote  the set of  labels of   old internal nodes (resp., singletons, young leaves) of $F$. 
 Let $S\subseteq \Oint(F)$ (resp.,~$S\subseteq \Yleaf(F)$)  and let $S_1, S_2, \ldots, S_m$ 	be the partition of $S$ such that $S_i  \subseteq  \Oint(T_i)$ (resp.,~$S_i\subseteq \Yleaf(T_i)$) for all $1\leq i\leq m$.
Relying on the transformation  $\Phi_x$, we  define the transformation   $\Phi_S$ by letting 
$$
\Phi_{S}(F)=(T'_1, T'_2, \ldots, T'_m) 
$$
where  $T'_i=\Phi_{S_i}(T_i)$ for all $i\in [m]$.
For example, let $F=(T_1, T_2, T_3)\in \mathcal{F}_9(3)$ in Fig.~\ref{PhiF} and let $S=\{3,5\}$. Then we have $\Phi_S(F)=(T'_1, T'_2, T'_3)$ as shown in Fig.~\ref{PhiF}, where $T'_1=\Phi_3(T_1)$, $T'_2=\Phi_5(T_2)$ and $T'_3=T_3$.

\begin{figure}[H]    
  	\begin{center}
  		\begin{tikzpicture}[font =\small , scale = 0.4]
  			\filldraw[fill=black](2*0.8,0.7*0)circle(0.1);
  			\filldraw[fill=black](0*0.8,0.7*-2)circle(0.1);
  			\filldraw[fill=black](2*0.8,0.7*-2)circle(0.1);
  			\filldraw[fill=black](4*0.8,0.7*-2)circle(0.1);
  			\filldraw[fill=black](0*0.8,0.7*-4)circle(0.1);
  			\filldraw[fill=black](2*0.8,0.7*-4)circle(0.1);
  			\draw (2*0.8,0.7*0)--(0*0.8,0.7*-2);
  			\draw (2*0.8,0.7*0)--(2*0.8,0.7*-2);
  			\draw (2*0.8,0.7*0)--(4*0.8,0.7*-2);
  			\draw (0*0.8,0.7*-4)--(0*0.8,0.7*-2);
  			\draw (2*0.8,0.7*-2)--(2*0.8,0.7*-4);
  			\coordinate [label={[text=black]above:$1$}] (x) at (2*0.8,0.7*0);
  			\coordinate [label={[text=black]below:$3$}] (x) at (2*0.8,0.7*-4);
  			\coordinate [label={[text=black]below:$2$}] (x) at (0*0.8,0.7*-4);
  			
  			\filldraw[fill=black](9*0.8,0.7*0)circle(0.1);
  			\filldraw[fill=black](6*0.8,0.7*-2)circle(0.1);
  			\filldraw[fill=black](9*0.8,0.7*-2)circle(0.1);
  			\filldraw[fill=black](12*0.8,0.7*-2)circle(0.1);
  			\filldraw[fill=black](8*0.8,0.7*-4)circle(0.1);
  			\filldraw[fill=black](10*0.8,0.7*-4)circle(0.1);
  			\filldraw[fill=black](13*0.8,0.7*-4)circle(0.1);
  			\filldraw[fill=black](11*0.8,0.7*-4)circle(0.1);
  			\draw (9*0.8,0.7*0)--(6*0.8,0.7*-2);
  			\draw (9*0.8,0.7*0)--(12*0.8,0.7*-2);
  			\draw (9*0.8,0.7*0)--(9*0.8,0.7*-2);
  			\draw (8*0.8,0.7*-4)--(9*0.8,0.7*-2);
  			\draw (10*0.8,0.7*-4)--(9*0.8,0.7*-2);
  			\draw (11*0.8,0.7*-4)--(12*0.8,0.7*-2);
  			\draw (13*0.8,0.7*-4)--(12*0.8,0.7*-2);
  			\coordinate [label={[text=black]above:$4$}] (x) at (9*0.8,0.7*0);
  			\coordinate [label={[text=black]below:$6$}] (x) at (8*0.8,0.7*-4);
  			\coordinate [label={[text=black]below:$8$}] (x) at (10*0.8,0.7*-4);
  			\coordinate [label={[text=black]below:$5$}] (x) at (11*0.8,0.7*-4);
  			\coordinate [label={[text=black]below:$7$}] (x) at (13*0.8,0.7*-4);
  			\coordinate [label=below:$ T_1$] (x) at (2*0.8,0.7*-6);
  			\coordinate [label=below:$ T_2$] (x) at (9*0.8,0.7*-6);
  			
  		\end{tikzpicture}
  	\end{center}
  	\caption{An example of a removable young leaf. }\label{rleaf}
  \end{figure}
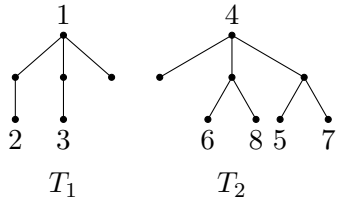

Given  a tree $T$, we simply write $u\in  T$ for any node $u$ of $T$.
\begin{definition}[Removable  young leaf]
\label{defiryleaf}
	Let $F=(T_1, T_2, \ldots, T_m)\in  \mathcal{F}_n(k)$.
	A young leaf $u\in T_m$ is said to be a  removable  young leaf of $F$ if  $u$ is a grand child of the root of $T_m$ and  the first $k-1$ children of the root  turn out to be leaves in   $\Phi_x(T_m)$, where $x$ is the label of $u$. 
\end{definition}
For instance,  the leaf with  label $5$ is a removable young leaf, whereas  the leaf with label $7$ is not a removable young leaf of  the forest $F=(T_1, T_2)$   in Fig.~\ref{rleaf}.

 Given a forest $F\in \mathcal{F}_{n}(k)$, let  $\rleaf(F)$ denote the total number of   removable old leaves and   removable young leaves of $F$. 
 Define
 $$ \overline{\mathcal{F}}^{*}_{n}(k)=\{F\mid F\in \overline{\mathcal{F}}_n(k),\,\,\ \yleaf(F)=0, \,\, 
 \rleaf(F)=0 \}  $$ 
 and
  $$ \widehat{\mathcal{F}}^{*}_{n}(k)=\{F\mid F\in \widehat{\mathcal{F}}_n(k),\,\,\ \yleaf(F)=0, \,\, 
\rleaf(F)=0 \}. $$ 

Given a forest $F=(T_1, T_2, \ldots, T_m)\in \mathcal{F}_{n}(k)$,  if   the root of $T_m$ has a grand child labeled by $x$ that is  an old internal node,  then set
  $\Oint^*(F)=\Oint(F)-\{x\}$; otherwise, set    $\Oint^*(F)=\Oint(F)$.  If $T_m$ is a singleton labeled by $y$, then define $\Si^*(F)=\Si(F)-\{y\}$; otherwise, we set $\Si^{*}(F)=\Si(F)$.

 We will record the properties  of $\Phi_{S}(F)$  in the following  lemmas.  
  \begin{lemma}\label{propPhiF0}
 	Given any forest  $F=(T_1, T_2, \ldots, T_m) \in \overline{\mathcal{F}}^{*}_{n}(k) $ and any subset $S \subseteq \Oint^*(F) $,  we have 	
 	$F'=\Phi_S(F)\in \overline{\mathcal{F}}_{n}(k) $ with
 	$$\rleaf(F')=0,\,\,\, \Si^*(F')=\Si^*(F)$$ and
 	$$\lleaf(F')-\si(F')=\lleaf(F)-\si(F)+|S|.$$ 
 \end{lemma}
 \pf
 Let $S_1, S_2, \ldots, S_m$	be the partition of $S$ such that $S_i  \subseteq  \Oint(T_i)$ for all $1\leq i\leq m$. Let $F'=(T'_1, T'_2, \ldots, T'_m)$ where $T'_i=\Phi_{S_i}(T_i)$. 
 By the definition of $\Phi_{S_i}(T_i)$,   one can easily check that $T'_i$ is a singleton if and only if $T_i$ is a singleton, and thus we have    $\Si(F')=\Si(F)$ and $\Si^*(F')=\Si^*(F)$. 
 If the root of $T_m$ has a grand child with label $x$ that is an old internal node, then we have  $x\notin S_m$ by the definition of $\Oint^*(F)$.  By the construction of $\Phi_{S_m}(T_m)$, we will do nothing for the root of $T_m$, the children of the root of $T_m$ and the grand children of the root of $T_m$ when we apply   $\Phi_{S_m}$ to $T_m$. This ensures that either $T'_m=T_m$ is a singleton or 
  the first $k-1$ children of the root of $T'_m=\Phi_{S_m}(T_m)$ are  leaves. 
  Hence, we have $F'\in \overline{\mathcal{F}}_{n}(k)$.   In view of Lemma \ref{lemPhi}, we have $ \Yleaf(F')=S$ and $\Oleaf(F')=\Oleaf(F)$.  Hence, we have 	$\lleaf(F')-\si(F')=\lleaf(F)-\si(F)+|S|$ since $\yleaf(F)=0$.  
 
 Now we proceed to show that $\rleaf(F')=0$.  Recall that  we will do nothing for the root of $T_m$, the children of the root of $T_m$ and the grand children of the root of $T_m$ when we apply   $\Phi_{S_m}$ to $T_m$. This ensures that there does not exist any removable old leaves or removable young leaves in $T'_m$ as $\rleaf(F)=0 $. 
 
 In order to show $\rleaf(F')=0$, it remains to show
   there does not exist any  removable old leaves in $T_i$ for all $i\in [m-1]$.  Take  $i\in [m-1]$ arbitrarily. 
Let $S'_i$ be the subset obtained from $S_i$ by removing the label of the old internal node that is a   grand child of the root of $T_i$ (if any).
 Clearly,   we do nothing for the root of $T_i$, the children of the root of $T_i$ and the grand children of the root of $T_i$ when we apply  $\Phi_{S'_i}$ to $T_i$.  Note that   $F$ does not contain any young leaves and any removable old leaves.  Hence, the tree $T''_i=\Phi_{S'_i}(T_i)$   does not contain  any removable old leaves. 
Then there does not exist any  removable old leaves in $T_i'=T_i''$ if $S'_i=S_i$. Otherwise,   
 suppose that we have   $\{x\}=S_i-S'_i$.  Note that the transformations $\Phi_y$  commute  for all $y\in S_i$. This yields that
 $T'_i=\Phi_{S_i}(T_i)=\Phi_{x}(\Phi_{S'_i}(T_i))=\Phi_{x}(T''_i)$. Then we will generate a young leaf   with label $x$  that is a grand child of  the root of $T'_i$ when we apply the  transformation $\Phi_{x}$ to $T''_i$.   
 Then,  by Definition~\ref{defioleaf},   there does not exist a  removable old leaf in $T'_i$.  Therefore, we have concluded that there does not exist any  removable old leaves in $T_i'$ for all $i\in [m-1]$,  competing the  proof. \qed

 \begin{lemma}\label{propPhiF}
 	Given any forest  $F=(T_1, T_2, \ldots, T_m) \in \widehat{\mathcal{F}}^{*}_{n}(k) $ and any $S \subseteq \Oint(F) $,  we have 	
 	$F'=\Phi_S(F)\in \widehat{\mathcal{F}}_{n}(k) $ with
 	$$\rleaf(F')=0,\,\,\, \Si^*(F')=\Si^*(F)$$ and
 	$$\lleaf(F')-\si(F')=\lleaf(F)-\si(F)+|S|.$$ 
 \end{lemma}
 \pf
 Let $S_1, S_2, \ldots, S_m$	be the partition of $S$ such that $S_i  \subseteq  \Oint(T_i)$ for all $1\leq i\leq m$. Let $F'=(T'_1, T'_2, \ldots, T'_m)$ where $T'_i=\Phi_{S_i}(T_i)$. 
 By the definition of $\Phi_{S_i}(T_i)$,   one can easily check that $T'_i$ is a singleton if and only if $T_i$ is a singleton, and thus we have    $\Si(F')=\Si(F)$ and $\Si^*(F')=\Si^*(F)$. Again by  the construction of $\Phi_{S_m}(T_m)$,    it is easily seen that the first $k-1$ children of the root of $T'_m$ are not all leaves. 
 This yields that $F'\in \widehat{\mathcal{F}}_{n}(k)$.    In view of Lemma~\ref{lemPhi}, we have $ \Yleaf(F')=S$ and $\Oleaf(F')=\Oleaf(F)$.  Hence, we have $\lleaf(F')-\si(F')=\lleaf(F)-\si(F)+|S|$ as $\yleaf(F)=0$. 
 
 Now we proceed to show that $\rleaf(F')=0$.   By the same reasoning as in the proof of Lemma \ref{propPhiF0}, one can   verify that there does not exist any removable old leaves in $T'_i$ for all $i\in [m-1]$.

 In order to show $\rleaf(F')=0$, it remains to show that there does not exist any removable young leaves or removable old leaves  in $T'_m$. Since $F'\in \widehat{\mathcal{F}}_{n}(k) $, the first $k-1$ children of the root of $T'_m$ are not all leaves. Then, by Definition \ref{defioleaf}, there does not exist any removable old leaf in $T'_m$. 
 
 Next we aim to show that there does not exist any removable young leaves   in $T'_m$. Let $S'_m$ be the subset obtained from $S_m$ by removing the label  of the old internal node that is a  grand child of the root of $T_m$ (if any).  Let $T''_m=\Phi_{S'_m}(T_m)$.  It is easily seen that 
 we will do nothing for the root of $T_m$, the children of the root of $T_m$  and the grand children of the root of $T_m$ when we apply  $\Phi_{S'_m}$ to $T_m$. This ensures that $T''_m$ does not  contain any removable young leaves since there does not exist any removable young leaves in $T_m$.

 If $S'_m=S_m$, then we have $T'_m=T''_m$ and thus we have 
  reached the conclusion that there does not exist any removable young leaves in $T'_m$. Otherwise, let $\{x\}=  S_m-S'_{m}$.
  Since $\Phi_x$ is an involution, we have  $\Phi_x(\Phi_x(T''_m))=T''_m$.    If  the node with label $x$ is a removable young leaf in $T'_m$, then the first $k-1$  children of the root of $\Phi_x(T'_m)=T''_m $ are all leaves  by Definition~\ref{defiryleaf}. However,  by the definition of $\Phi_{S'_m}(T_m)$, the first $k-1$ children  of the root of $T''_m$ are not all leaves as $F\in\widehat{F}_n(k)$, a contradiction. 
  Hence, the node with label $x$ is not a removable young leaves in $T'_m$.   As $\yleaf(F)=0$,   $x$ is the smallest label among all  the labels of  young leaves  that are grand children of the root of $T'_m$.  Suppose that the node with label $y$ is a removable young leaf in $T'_m$ with $y\neq x$. If the leaf with label  $x$ lies in  the first $k-1$ subtrees of the root of $T'_m$, then one can easily check that  the first $k-1$  children of the root of  $\Phi_y(T'_m)$ are not all leaves and thus the leaf with label $y$ is not a   removable young leaf.
 If the leaf $x$ lies in the rightmost subtree of the root of $T'_m$,  then the leaf with label $y$ is  a removable young leaf  would imply that  the leaf with label $x$ is   a removable young leaf since $y>x$. Therefore, we have deduced that  $T'_m$ does not contains any removable young leaves, completing  the   proof. \qed

 Define
  $$
  \mathcal{X}_{n}(k)=\{(F,S)\mid F\in \mathcal{F}_n(k), \,\,\yleaf(F)=0,   \,\,\, S\subseteq \Oint(F)\bigcup \Si^*(F)  \}, 
 $$
 and 
  $$
 \mathcal{Y}_{n}(k)=\{(F,S)\mid F\in \mathcal{F}_n(k), \,\,   S\subseteq  \Si^*(F)   \}.
 $$

\noindent{\bf The map $\Theta$ from $\mathcal{X}_n(k)$ to $\mathcal{Y}_n(k)$.}\\
 Given an  ordered pair $(F,S)\in  \mathcal{X}_{n}(k) $, let  $(S_1, S_2)$ be the partition of  $S$  such that   $S_1\subseteq  \Oint(F)$ and  $S_2  \subseteq  \Si^*(F)$.  Set $\Theta(F,S)=(\Phi_{S_1}(F), S_2)$.  
 
  It is apparent that we have $F'=\Phi_{S_1}(F)\in \mathcal{F}_{n}(k)$ and     $\Si^*(F')=\Si^*(F)$. This ensures that $S_2\subseteq \Si^*(F')$ and thus we have $\Theta(F,S)\in \mathcal{Y}_n(k)$. 

\noindent{\bf The map $\Theta'$ from $\mathcal{Y}_n(k)$ to $\mathcal{X}_n(k)$.}\\
Given an ordered pair $(F, S_1)\in \mathcal{Y}_{n}(k) $, let $S_2$ be the set of the labels of young leaves of 
$F$. Set $\Theta'(F, S_1)=(\Phi_{S_2}(F), S_1\bigcup S_2)$.  

By Lemma~\ref{lemPhi}, the resulting forest $F'=\Phi_{S_2}(F)$ does not contain any young leaves. Moreover, we have $S_1\subseteq \Si^*(F')$ and $S_2\subseteq \Oint(F')$.  Therefore, we have  $\Theta'(F, S_1)\in \mathcal{X}_n(k) $, and thus the map $\Theta'$ is well defined.

  \begin{proposition} \label{Thetageneral}
 	For $n,k\geq 1$, the maps $\Theta$ and $\Theta'$ induce a one-to-one correspondence between $\mathcal{X}_n(k)$ and $\mathcal{Y}_n(k)$. 
 \end{proposition} 
 \pf  It remains to  show that the maps $\Theta$ and $\Theta'$ are inverses of each other. Here we only deal with the equality   $\Theta'(\Theta(F,S))=(F,S)$ for any  $(F,S)\in \mathcal{X}_{n}(k) $. By similar arguments, one can verify that the equality $\Theta(\Theta'(F,S))=(F,S)$ holds for any  $(F,S)\in \mathcal{Y}_n(k) $.   Recall that
 $\Theta(F,S)=(\Phi_{S_1}(F), S_2)$, where  $(S_1, S_2)$ is the partition of  $S$  such that   $S_1\subseteq  \Oint(F)$ and  $S_2  \subseteq  \Si^*(F)$.   
 In view of Lemma~\ref{lemPhi}, one can easily check that    $\Yleaf(\Phi_{S_1}(F))=S_1$.  Then we have  $\Theta'(\Phi_{S_1}(F), S_2)=(F', S')$ where $F'=\Phi_{S_1}(\Phi_{S_1}(F))$ and $S'=S_1\cup S_2$.  Since the transformations $\Phi_x$ are involutions  and commute for any $x\in [n]$, we have $F'=\Phi_{S_1}(\Phi_{S_1}(F))=F$. Hence, we have concluded that the equality
 $\Theta'(\Theta(F,S))=(F,S)$ holds for any  $(F,S)\in  \mathcal{X}_{n}(k) $, completing the proof. \qed

Define 
 $$
 \overline{\mathcal{X}}_{n}(k)=\{(F,S)\mid F\in \overline{\mathcal{F}}^*_n(k), \,\,  S\subseteq \Oint^*(F)\bigcup \Si^*(F)  \}, 
 $$
   $$
 \overline{\mathcal{Y}}_{n}(k)=\{(F,S)\mid F\in \overline{\mathcal{F}}_n(k),\,\, \rleaf(F)=0,  \,\, \,\, S\subseteq  \Si^*(F)  \}, $$
  $$
 \widehat{\mathcal{X}}_{n}(k)=\{(F,S)\mid F\in \widehat{\mathcal{F}}^*_n(k), \,\,   S\subseteq \Oint(F)\bigcup \Si^*(F)  \}, 
 $$
 and 
 $$
 \widehat{\mathcal{Y}}_{n}(k)=\{(F,S)\mid F\in \widehat{\mathcal{F}}_n(k),\,\, \rleaf(F)=0,  \,\, \,\, S\subseteq  \Si^*(F)  \}. 
 $$
 Clearly, we have $\overline{\mathcal{X}}_{n}(k),  \widehat{\mathcal{X}}_{n}(k)\subseteq\mathcal{X}_{n}(k)$ and $\overline{\mathcal{Y}}_{n}(k),  \widehat{\mathcal{Y}}_{n}(k)\subseteq\mathcal{Y}_{n}(k)$.  
 
   \begin{proposition} \label{Theta0}
  	For $n,k\geq 1$, the maps $\Theta$ and $\Theta'$ induce a one-to-one correspondence between $\overline{\mathcal{X}}_{n}(k)$ and $\overline{\mathcal{Y}}_{n}(k)$ so that  for any  $(F,S)\in \overline{\mathcal{X}}_{n}(k) $, we have 
  	$\Theta(F,S)=(\Phi_{S_1}(F), S_2)$ with  	
  	$$
  	\lleaf(\Phi_{S_1}(F))-\si(\Phi_{S_1}(F))=\lleaf(F)-\si(F)+|S_1|,
  	$$
  	where   $(S_1, S_2)$ is the partition of  $S$  such that   $S_1\subseteq  \Oint^*(F)$ and  $S_2  \subseteq  \Si^*(F)$.  
  \end{proposition} 
  \pf
  By Lemma \ref{propPhiF0}, we have $\Phi_{S_1}(F)\in \overline{\mathcal{F}}_n(k)$ with  
  $$\rleaf(\Phi_{S_1}(F))=0,\,\,\, \Si^*(\Phi_{S_1}(F) )=\Si^*(F)$$ and
  $$\lleaf(\Phi_{S_1}(F))-\si(\Phi_{S_1}(F))=\lleaf(F)-\si(F)+|S_1|.$$  
  Hence, we have $\Theta(F,S)\in \overline{\mathcal{Y}}_{n,k}$ and thus the map  $\Theta$ is   well-defined.

  In the following, we aim to show that $\Theta'$ is a map from $\overline{\mathcal{Y}}_{n}(k) $ to $\overline{\mathcal{X}}_{n}(k) $. 
  Given an ordered pair $(F, S_1)\in\overline{\mathcal{Y}}_{n}(k) $, let $S_2$ be the set of the labels of young leaves of 
  $F$. We have  $\Theta'(F, S_1)=(F',S')=(\Phi_{S_2}(F), S_1\bigcup S_2)$. According to the construction of $\Phi_{S_2}(F)$, one can easily check that $F'\in \overline{\mathcal{F}}_n(k)$.   By Proposition~\ref{Thetageneral},   
  in order to show that $\Theta'(F, S_1)\in \overline{\mathcal{X}}_{n}(k)$, 
  it remains to show that $S_2\subseteq \Oint^*(F')$ and  $\rleaf(F')=0$.
  
  First we aim to show that $S_2\subseteq \Oint^*(F')$.
  Let $F=(T_1, T_2, \ldots, T_m)$. Since $F\in \overline{\mathcal{F}}_n(k)$ and $\rleaf(F)=0$, 
  the root of $T_m$ does not have a grand child that is a young leaf by Definition \ref{defiryleaf}. This ensures that $S_2\subseteq \Oint^*(F')$  by the construction of $\Phi_{S_2}(F)$.

  Next we proceed to show    that  $\rleaf(F')=0$.
  Since all the transformations $\Phi_x$ are involutions and commute, we have
  $\Phi_{S_2}(F')=\Phi_{S_2}(\Phi_{S_2}(F))=F$. 
  Let $F'=(T'_1, T'_2, \ldots, T'_m)$. Assume that $u\in T'_i$ is a removable old leaf of $F'$. This ensures that  all the labels of the grand children of the  root of $T'_i$ do not belong to the set $S_2$ since $S_2\subseteq \Oint(F') $ by Lemma \ref{lemPhi}. Then when we apply the transformation $\Phi_{S_2}$ to $F'$,  the node $u$ will  also be a removable old leaf in the resulting forest $F$. This contradicts the fact that $\rleaf(F)=0$. Therefore, we have  $\rleaf(F')=0$, and thus the map $\Theta'$ is well defined. 
  
  By Proposition~\ref{Thetageneral}, the maps $\Theta$ and $\Theta'$ are inverses of each other. Therefore, we have concluded that the maps $\Theta$ and $\Theta'$ induce a bijection  between $\overline{\mathcal{X}}_n(k)$ and $\overline{\mathcal{Y}}_n(k)$, completing the proof. 
  \qed

  \begin{proposition} \label{Theta}
  	For $n,k\geq 1$, the maps $\Theta$ and $\Theta'$ induce a one-to-one correspondence between $\widehat{\mathcal{X}}_{n}(k)$ and $\widehat{\mathcal{Y}}_{n}(k)$ so that  for any  $(F,S)\in \widehat{\mathcal{X}}_{n}(k) $, we have 
  	$\Theta(F,S)=(\Phi_{S_1}(F), S_2)$ with  	
  	$$
  	\lleaf(\Phi_{S_1}(F))-\si(\Phi_{S_1}(F))=\lleaf(F)-\si(F)+|S_1|, 
  	$$
  	where   $(S_1, S_2)$ is the partition of  $S$  such that   $S_1\subseteq  \Oint(F)$ and  $S_2  \subseteq  \Si^*(F)$.  
  \end{proposition}
 \pf  By Lemma~\ref{propPhiF}, we have $\Phi_{S_1}(F)\in \widehat{\mathcal{F}}_n(k)$ with  
 $$\rleaf(\Phi_{S_1}(F))=0,\,\,\, \Si^*(\Phi_{S_1}(F) )=\Si^*(F)$$ and
 $$\lleaf(\Phi_{S_1}(F))-\si(\Phi_{S_1}(F))=\lleaf(F)-\si(F)+|S_1|.$$   
 Hence, we have $\Theta(F,S)\in \widehat{\mathcal{Y}}_{n,k}$ and thus the map  $\Theta$ is   well-defined. 
 
  In the following, we aim to show that $\Theta'$ is a map from $\widehat{\mathcal{Y}}_{n}(k)$ to $\widehat{\mathcal{X}}_{n}(k)$, which will complete the proof in view of Proposition~\ref{Thetageneral}. 
 Given an ordered pair $(F, S_1)\in\widehat{\mathcal{Y}}_{n}(k) $, let $S_2$ be the set of the labels of young leaves of 
 $F$. We have  $\Theta'(F, S_1)=(F',S')=(\Phi_{S_2}(F), S_1\bigcup S_2)$. Since $F$ does not contain any removable young leaves, we have $\Phi_{S_2}(F)\in \widehat{\mathcal{F}}_n(k)$.  By Proposition~\ref{Thetageneral},   
 in order to show that $\Theta'(F, S_1)\in \widehat{\mathcal{X}}_{n}(k)$, it remains to show that $\rleaf(F')=0$. This can be justified by similar arguments as in the proof of Proposition~\ref{Theta0} and the proof is completed.  
   \qed

 \subsection{Introducing two transformations on $\mathcal{F}_n(k)$}

Let $F=(T_1, T_2, \ldots, T_m)\in \mathcal{F}_n(k)$. 
 Assume that the node $u\in T_i$ is labeled by $x$ in $F$. 
 We define the fundamental transformation
$\Psi_x(F)$ as follows:
\begin{itemize}
	\item If   $u$ is a singleton and $u\notin T_m$, choose the least integer $j>i$ such that $T_j$ is either  a singleton or $T_j=T_m$. 
	 
	\begin{itemize}
	\item[(i)] If $T_j$ is a singleton, then attach $k$ unlabeled leaves to $u$, say $u_1, u_2, \ldots, u_k$.  Attach  the trees between $T_i$ and $T_j$ (including $T_j$)  to the node $u_k$ as subtrees from left to right. Denote by $\Psi_x(F)$ the  resulting forest.  For instance, let $F_1$ and $F_2$ be the forests illustated  in Fig.~\ref{Psi_2}. Then we have $\Psi_2(F_1)=F_2$. 
	\begin{figure}
	   \begin{center}	\begin{tikzpicture}[font =\small , scale = 0.3]        \filldraw[fill=black](2,0)circle(0.1);        \filldraw[fill=black](0,-2)circle(0.1);        \filldraw[fill=black](2,-2)circle(0.1);        \filldraw[fill=black](4,-2)circle(0.1);        \filldraw[fill=black](0,-4)circle(0.1);        \filldraw[fill=black](4,-4)circle(0.1);        \filldraw[fill=red](10,-2)circle(0.1);        \coordinate [label=below:$3$] (x) at (0,-4);        \coordinate [label=below:$7$] (x) at (4,-4);        \draw(0,-2)--(2,0);        \draw(2,-2)--(2,0);        \draw(4,-2)--(2,0);        \draw(0,-4)--(0,-2);        \draw(4,-4)--(4,-2);        \draw[red](10,0)--(10,-2);        \filldraw[fill=red](6,0)circle(0.1);        \coordinate [label={[text=red]above:$2$}] (x) at (6,0);         \filldraw[fill=red](10,0)circle(0.1);         \coordinate [label={[text=red]above:$4$}] (x) at (10,0);         \coordinate [label={[text=red]below:$6$}] (x) at (12,-4);         \filldraw[fill=red](8,-2)circle(0.1);        \filldraw[fill=red](12,-2)circle(0.1);        \filldraw[fill=red](12,-2)circle(0.1);        \filldraw[fill=red](12,-4)circle(0.1);        \draw[red] (8,-2)--(10,0);        \draw[red] (12,-2)--(10,0);        \draw[red] (12,-4)--(12,-2);        \filldraw[fill=red](14,0)circle(0.1);        \coordinate [label={[text=red]above:$5$}] (x) at (14,0);        \filldraw(16,0)circle(0.1);        \coordinate [label={[text=black]above:$8$}] (x) at (16,0);        \coordinate [label=below:$ F_1$] (x) at (10,-5);        \draw[line width=1.2pt][->] (9,-6) -- (9,-8);        \filldraw[fill=black](2,-9)circle(0.1);        \filldraw[fill=black](0,-11)circle(0.1);        \filldraw[fill=black](2,-11)circle(0.1);        \filldraw[fill=black](4,-11)circle(0.1);        \filldraw[fill=black](0,-13)circle(0.1);        \filldraw[fill=black](4,-13)circle(0.1);        \coordinate [label=below:$3$] (x) at (0,-13);        \coordinate [label=below:$7$] (x) at (4,-13);        \draw(0,-11)--(2,-9);        \draw(2,-11)--(2,-9);        \draw(4,-11)--(2,-9);        \draw(0,-13)--(0,-11);        \draw(4,-13)--(4,-11);        \filldraw[fill=red](10,-9)circle(0.1);         \coordinate [label={[text=red]above:$2$}] (x) at (10,-9);         \filldraw[fill=red](8,-11)circle(0.1);         \filldraw[fill=red](10,-11)circle(0.1);         \filldraw[fill=red](12,-11)circle(0.1);         \filldraw[fill=red](10,-13)circle(0.1);         \filldraw[fill=red](14,-13)circle(0.1);         \filldraw[fill=red](8,-15)circle(0.1);         \filldraw[fill=red](10,-15)circle(0.1);         \filldraw[fill=red](12,-15)circle(0.1);         \filldraw[fill=red](12,-17)circle(0.1);        \draw[red] (10,-9)--(10,-11);        \draw[red] (10,-9)--(8,-11);        \draw[red] (10,-9)--(12,-11);        \draw[red] (10,-13)--(12,-11);        \draw[red] (14,-13)--(12,-11);        \draw[red] (12,-15)--(10,-13);        \draw[red] (8,-15)--(10,-13);        \draw[red] (10,-15)--(10,-13);        \draw[red] (12,-15)--(12,-17);        \filldraw[fill=black](16,-9)circle(0.1);        \coordinate [label={[text=red]below:$5$}] (x) at (14,-13);        \coordinate [label={[text=red]above:$4$}] (x) at (10,-13);        \coordinate [label={[text=red]below:$6$}] (x) at (12,-17);        \coordinate [label={[text=black]above:$8$}] (x) at (16,-9);        \coordinate [label=below:$ F_2$] (x) at (10,-18);                \coordinate [label=above:$1$] (x) at (2,-9);        \coordinate [label=above:$1$] (x) at (2,0);                \end{tikzpicture}    \end{center}
		\caption{The   transformation  $\Psi_{2}$ applying to $F_1$. }\label{Psi_2}
	\end{figure}
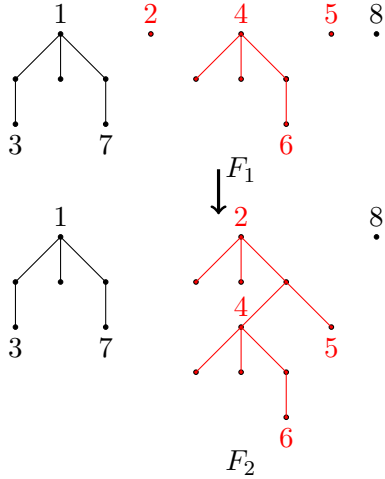
	\item[(ii)] Otherwise,  we must have $T_j=T_m$.   Suppose that   $T_m$ is rooted at  the node $v$ with label $y$ and   that the node  $v'$  is the rightmost child of $v$. 
	 Attach  the trees between $T_i$ and $T_m$  to the node $v'$ as   subtrees  and  relabel   $v$ by $x$. Then attach a leaf labeled by $y$ to   $v'$ and rearrange the subtrees  of  $v'$ so that the children of $v_k$ are increasing from left to right.	 Denote by $\Psi_x(F)$ the  resulting forest.   For instance, let $F_3$ and $F_4$ be the forests illustated  in Fig.~\ref{Psi_22}. Then we have $\Psi_2(F_3)=F_4$.
	 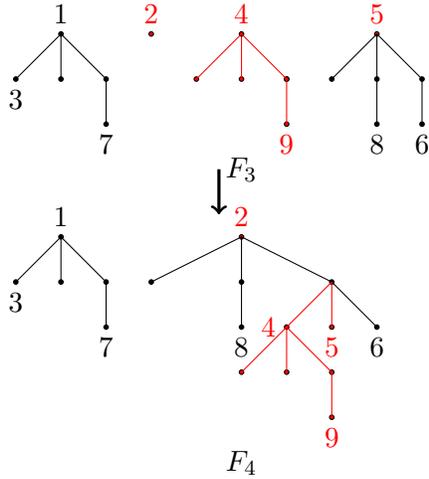
\begin{figure}    \begin{center}	\begin{tikzpicture}[font =\small , scale = 0.3]        \filldraw[fill=black](2,0)circle(0.1);        \filldraw[fill=black](0,-2)circle(0.1);        \filldraw[fill=black](2,-2)circle(0.1);        \filldraw[fill=black](4,-2)circle(0.1);        \filldraw[fill=black](4,-4)circle(0.1);        \filldraw[fill=red](10,-2)circle(0.1);        \coordinate [label=below:$3$] (x) at (0,-2);        \coordinate [label=below:$7$] (x) at (4,-4);        \coordinate [label=above:$1$] (x) at (2,0);        \draw(0,-2)--(2,0);        \draw(2,-2)--(2,0);        \draw(4,-2)--(2,0);        \draw(4,-4)--(4,-2);        \draw[red](10,0)--(10,-2);        \filldraw[fill=red](6,0)circle(0.1);        \coordinate [label={[text=red]above:$2$}] (x) at (6,0);         \filldraw[fill=red](10,0)circle(0.1);         \coordinate [label={[text=red]above:$4$}] (x) at (10,0);         \coordinate [label={[text=red]below:$9$}] (x) at (12,-4);         \filldraw[fill=red](8,-2)circle(0.1);        \filldraw[fill=red](12,-2)circle(0.1);        \filldraw[fill=red](12,-2)circle(0.1);        \filldraw[fill=red](12,-4)circle(0.1);        \draw[red] (8,-2)--(10,0);        \draw[red] (12,-2)--(10,0);        \draw[red] (12,-4)--(12,-2);        \filldraw[fill=red](16,0)circle(0.1);        \filldraw[fill=black](16,-2)circle(0.1);        \filldraw[fill=black](18,-2)circle(0.1);        \filldraw[fill=black](14,-2)circle(0.1);        \filldraw[fill=black](16,-4)circle(0.1);        \filldraw[fill=black](18,-4)circle(0.1);        \coordinate [label=below:$8$] (x) at (16,-4);        \coordinate [label=below:$6$] (x) at (18,-4);        \coordinate [label={[text=red]above:$5$}] (x) at (16,0);        \draw(16,-2)--(16,0);        \draw(18,-2)--(16,0);        \draw(14,-2)--(16,0);        \draw(16,-4)--(16,-2);        \draw(18,-4)--(18,-2);                \coordinate [label=below:$ F_3$] (x) at (10,-5);        \draw[line width=1.2pt][->] (9,-6) -- (9,-8);        \filldraw[fill=black](2,-9)circle(0.1);        \filldraw[fill=black](0,-11)circle(0.1);        \filldraw[fill=black](2,-11)circle(0.1);        \filldraw[fill=black](4,-11)circle(0.1);        \filldraw[fill=black](4,-13)circle(0.1);        \coordinate [label=below:$3$] (x) at (0,-11);        \coordinate [label=above:$1$] (x) at (2,-9);        \coordinate [label=below:$7$] (x) at (4,-13);        \draw(0,-11)--(2,-9);        \draw(2,-11)--(2,-9);        \draw(4,-11)--(2,-9);        \draw(4,-13)--(4,-11);        \filldraw[fill=red](10,-9)circle(0.1);         \coordinate [label={[text=red]above:$2$}] (x) at (10,-9);         \filldraw[fill=black](6,-11)circle(0.1);         \filldraw[fill=black](10,-11)circle(0.1);         \filldraw[fill=black](14,-11)circle(0.1);         \filldraw[fill=black](10,-13)circle(0.1);         \filldraw[fill=black](16,-13)circle(0.1);        \draw(10,-9)--(10,-11);        \draw(6,-11)--(10,-9);        \draw(14,-11)--(10,-9);        \draw(16,-13)--(14,-11);        \draw(10,-13)--(10,-11);        \coordinate [label=below:$8$] (x) at (10,-13);        \coordinate [label=below:$6$] (x) at (16,-13);        \filldraw[fill=red](12,-13)circle(0.1);        \filldraw[fill=red](14,-13)circle(0.1);        \filldraw[fill=red](12,-15)circle(0.1);        \filldraw[fill=red](10,-15)circle(0.1);        \filldraw[fill=red](14,-15)circle(0.1);        \filldraw[fill=red](14,-17)circle(0.1);        \draw[red] (14,-11)--(14,-13);        \draw[red] (14,-11)--(12,-13);        \draw[red] (12,-15)--(12,-13);        \draw[red] (14,-15)--(12,-13);        \draw[red] (10,-15)--(12,-13);        \draw[red] (14,-17)--(14,-15);        \coordinate [label={[text=red]left:$4$}] (x) at (12,-13);        \coordinate [label={[text=red]below:$5$}] (x) at (14,-13);        \coordinate [label={[text=red]below:$9$}] (x) at (14,-17);        \coordinate [label=below:$ F_4$] (x) at (10,-18);        \end{tikzpicture}    \end{center}
	 		\caption{The   transformation $\Psi_{2}$ applying to $F_3$. }\label{Psi_22}
	 \end{figure}
 
	 	\end{itemize}
	\item If   $u$ is a removable old  leaf, then suppose that $T_i$ is rooted at the node $v$ with label $y$ and that the rightmost child of $v$ is the node $v'$.  Remove all the subtrees of $v'$ 
	  and put  them  immediately to the right of $T_i$ from left to right in the same order. Then remove all the leaves of $v$  together with the edges incident to them.    Denote by $\Psi_x(F)$ the  resulting forest.  For instance,  if we let $F_1$ and $F_2$ be the forests illustrated  in Fig.~\ref{Psi_2}, then we have  $\Psi_5(F_2)=F_1$. 
	
	\item If  $u$ is a removable young leaf,  then we have $u\in T_m$. Suppose that   $T_m$ is rooted at  the node $v$ with label $y$ and   that the node  $v'$  is the rightmost child of $v$.
	 Remove all  subtrees, say $T'_1, T'_2, \ldots, T'_\ell$ from left to right,  of $v'$ that locate to the left of  the edge between $v'$ and $u$. Put $T'_1, T'_2, \ldots, T'_\ell$  immediately to the left  of $T_m$ from left to right.  Then remove   $u$ together with the edge incident to it and relabel the root $v$ by $x$. Insert a singleton with label $y$ immediately before   $T'_1$.   Denote by $\Psi_x(F)$ the  resulting forest.   For instance,  if we let $F_3$ and $F_4$ be the forests illustrated  in Fig.~\ref{Psi_22}, then we have  $\Psi_5(F_4)=F_3$. 
	
	\item  Otherwise, set $\Psi_x(F)=F$. 
\end{itemize}

In the following, we record the basic properties of this transformation  $\Psi_x(F)$,  which follow immediately from its  definition.

\begin{lemma}\label{lemPsi}
	Let $F=(T_1, T_2, \ldots, T_m)\in \mathcal{F}_n(k)$. 
	Assume that the node $u\in T_i$ is labeled by $x$ in $F$.  Then we have $\Psi_x(F)\in \mathcal{F}_n(k) $  and $F\in \overline{\mathcal{F}}_n(k)$ if and only  if $\Psi_x(F)\in \overline{\mathcal{F}}_n(k)$.
	Moreover, we have 
 	\begin{equation}\label{eqPsi2}
		\lleaf(\Psi_x(F))-\si(\Psi_x(F))=	\lleaf(F)-\si(F)+1
		\end{equation}
	when $u$ is a singleton and $u\notin T_m$, and we have
		\begin{equation}\label{eqPsi3}
		\lleaf(\Psi_x(F))-\si(\Psi_x(F))=	\lleaf(F)-\si(F)-1
	\end{equation}
	when  $u$ is either a removable old leaf or a removable young leaf. 
\end{lemma}

\noindent{\bf The tansformation $\alpha$. }\\
Given a pair $(F,S)  $  with $S\subseteq \Si^*(F)$, we define a  transformation $\alpha(F,S)$  as follows.  Let $x$ be the \red{\bf greatest} element of the set $S$.  Set $\alpha(F,S)=(F', S')$, where $F'=\Psi_x(F)$ and $S'=S-\{x\}$. 

Let $i$ be a nonnegative integer. Denote by 
$\alpha^i(F, S)$ the resulting pair  obtained after the $i$-th application of $\alpha$ to  the pair $(F,S)$. Clearly, we have $\alpha^0(F,S)=(F,S)$.

\noindent{\bf The map $\Gamma$ from   $\mathcal{Y}_n(k)$  to $\mathcal{F}_n(k)$. }\\
 Given an ordered pair $(F,S)\in \mathcal{Y}_{n}(k)$, let $S=\{i_1, i_2, \ldots, i_\ell\}  \subseteq  \Si^*(F)$ with $i_1>i_2>\cdots>i_{\ell}$.  After the $\ell$-th application of $\alpha$ to $(F,S)$, we will reach a pair $(F^{(\ell)}, S^{(\ell)})$ with $S^{(\ell)}=\emptyset.$  Set $\Gamma(F,S)=F^{(\ell)}$.

\begin{figure}    \begin{center}	\begin{tikzpicture}[font =\small , scale = 0.4]        \filldraw[fill=black](0*0.8,0.7*1)circle(0.1);        \filldraw[fill=black](4*0.8,0.7*1)circle(0.1);        \filldraw[fill=black](2*0.8,0.7*-1)circle(0.1);        \filldraw[fill=black](4*0.8,0.7*-1)circle(0.1);        \filldraw[fill=black](6*0.8,0.7*-1)circle(0.1);        \filldraw[fill=black](4*0.8,0.7*-3)circle(0.1);        \draw(4*0.8,0.7*1)--(4*0.8,0.7*-1);        \draw(4*0.8,0.7*-1)--(4*0.8,0.7*-3);        \draw(2*0.8,0.7*-1)--(4*0.8,0.7*1);        \draw(6*0.8,0.7*-1)--(4*0.8,0.7*1);        \coordinate [label=above:$1$] (x) at (0*0.8,0.7*1);        \coordinate [label=above:$2$] (x) at (4*0.8,0.7*1);        \coordinate [label=below:$5$] (x) at (4*0.8,0.7*-3);                \filldraw[fill=red](8*0.8,0.7*1)circle(0.1);        \filldraw[fill=red](16*0.8,0.7*1)circle(0.1);        \filldraw[fill=red](12*0.8,0.7*1)circle(0.1);        \filldraw[fill=red](12*0.8,0.7*-1)circle(0.1);        \filldraw[fill=red](10*0.8,0.7*-1)circle(0.1);        \filldraw[fill=red](14*0.8,0.7*-1)circle(0.1);        \filldraw[fill=red](14*0.8,0.7*-3)circle(0.1);        \draw[red] (12*0.8,0.7*1)--(12*0.8,0.7*-1);        \draw[red] (12*0.8,0.7*1)--(10*0.8,0.7*-1);        \draw[red] (12*0.8,0.7*1)--(14*0.8,0.7*-1);        \draw[red] (14*0.8,0.7*-3)--(14*0.8,0.7*-1);`       \coordinate [label={[text=red]above:$3$}] (x) at (8*0.8,0.7*1);        \coordinate [label={[text=red]above:$4$}] (x) at (12*0.8,0.7*1);        \coordinate [label={[text=red]below:$7$}] (x) at (14*0.8,0.7*-3);        \coordinate [label={[text=red]above:$6$}] (x) at (16*0.8,0.7*1);        \filldraw[fill=black](20*0.8,0.7*1)circle(0.1);        \filldraw[fill=black](20*0.8,0.7*-1)circle(0.1);        \filldraw[fill=black](18*0.8,0.7*-1)circle(0.1);        \filldraw[fill=black](22*0.8,0.7*-1)circle(0.1);        \filldraw[fill=black](18*0.8,0.7*-3)circle(0.1);        \filldraw[fill=black](22*0.8,0.7*-3)circle(0.1);        \draw(20*0.8,0.7*1)--(20*0.8,0.7*-1);        \draw(20*0.8,0.7*1)--(18*0.8,0.7*-1);        \draw(20*0.8,0.7*1)--(22*0.8,0.7*-1);        \draw(18*0.8,0.7*-3)--(20*0.8,0.7*-1);        \draw(22*0.8,0.7*-3)--(20*0.8,0.7*-1);        \coordinate [label=above:$8$] (x) at (20*0.8,0.7*1);        \coordinate [label=below:$9$] (x) at (18*0.8,0.7*-3);        \coordinate [label=below:$10$] (x) at (22*0.8,0.7*-3);                \draw[line width=1.2pt][->] (11*0.8,0.7*-6) -- (11*0.8,0.7*-8);      
  			\filldraw[fill=black](0*0.8,0.7*-10)circle(0.1);        \filldraw[fill=black](4*0.8,0.7*-10)circle(0.1);        \filldraw[fill=black](2*0.8,0.7*-12)circle(0.1);        \filldraw[fill=black](4*0.8,0.7*-12)circle(0.1);        \filldraw[fill=black](6*0.8,0.7*-12)circle(0.1);        \filldraw[fill=black](4*0.8,0.7*-14)circle(0.1);        \draw(4*0.8,0.7*-10)--(4*0.8,0.7*-12);        \draw(4*0.8,0.7*-12)--(4*0.8,0.7*-14);        \draw(2*0.8,0.7*-12)--(4*0.8,0.7*-10);        \draw(6*0.8,0.7*-12)--(4*0.8,0.7*-10);        \coordinate [label=above:$1$] (x) at (0*0.8,0.7*-10);        \coordinate [label=above:$2$] (x) at (4*0.8,0.7*-10);        \coordinate [label=below:$5$] (x) at (4*0.8,0.7*-14);        \filldraw[fill=red](11*0.8,0.7*-10)circle(0.1);        \filldraw[fill=red](11*0.8,0.7*-12)circle(0.1);        \filldraw[fill=red](13*0.8,0.7*-12)circle(0.1);        \filldraw[fill=red](9*0.8,0.7*-12)circle(0.1);        \filldraw[fill=red](14*0.8,0.7*-14)circle(0.1);        \filldraw[fill=red](11*0.8,0.7*-14)circle(0.1);        \filldraw[fill=red](11*0.8,0.7*-16)circle(0.1);        \filldraw[fill=red](9*0.8,0.7*-16)circle(0.1);        \filldraw[fill=red](13*0.8,0.7*-16)circle(0.1);        \filldraw[fill=red](13*0.8,0.7*-18)circle(0.1);        \draw[red] (11*0.8,0.7*-10)--(11*0.8,0.7*-12);        \draw[red] (11*0.8,0.7*-10)--(9*0.8,0.7*-12);        \draw[red] (11*0.8,0.7*-10)--(13*0.8,0.7*-12);        \draw[red] (14*0.8,0.7*-14)--(13*0.8,0.7*-12);        \draw[red] (11*0.8,0.7*-14)--(13*0.8,0.7*-12);        \draw[red] (11*0.8,0.7*-14)--(9*0.8,0.7*-16);        \draw[red] (11*0.8,0.7*-14)--(11*0.8,0.7*-16);        \draw[red] (11*0.8,0.7*-14)--(13*0.8,0.7*-16);        \draw[red] (13*0.8,0.7*-18)--(13*0.8,0.7*-16);        \coordinate [label={[text=red]above:$3$}] (x) at (11*0.8,0.7*-10);        \coordinate [label={[text=red]left:$4$}] (x) at (11*0.8,0.7*-14);        \coordinate [label={[text=red]below:$6$}] (x) at (14*0.8,0.7*-14);        \coordinate [label={[text=red]below:$7$}] (x) at (13*0.8,0.7*-18);        \filldraw[fill=black](20*0.8,0.7*-10)circle(0.1);        \filldraw[fill=black](20*0.8,0.7*-12)circle(0.1);        \filldraw[fill=black](18*0.8,0.7*-12)circle(0.1);        \filldraw[fill=black](22*0.8,0.7*-12)circle(0.1);        \filldraw[fill=black](18*0.8,0.7*-14)circle(0.1);        \filldraw[fill=black](22*0.8,0.7*-14)circle(0.1);        \draw(20*0.8,0.7*-10)--(20*0.8,0.7*-12);        \draw(20*0.8,0.7*-10)--(18*0.8,0.7*-12);        \draw(20*0.8,0.7*-10)--(22*0.8,0.7*-12);        \draw(18*0.8,0.7*-14)--(20*0.8,0.7*-12);        \draw(22*0.8,0.7*-14)--(20*0.8,0.7*-12);        \coordinate [label=above:$8$] (x) at (20*0.8,0.7*-10);        \coordinate [label=below:$9$] (x) at (18*0.8,0.7*-14);        \coordinate [label=below:$10$] (x) at (22*0.8,0.7*-14);     
  			\draw[line width=1.2pt][->] (11*0.8,0.7*-22) -- (11*0.8,0.7*-24); 
  			\filldraw[fill=red](11*0.8,0.7*-26)circle(0.1);        \filldraw[fill=black](5*0.8,0.7*-28)circle(0.1);        \filldraw[fill=red](17*0.8,0.7*-28)circle(0.1);        \filldraw[fill=black](11*0.8,0.7*-28)circle(0.1);        \filldraw[fill=black](9*0.8,0.7*-30)circle(0.1);        \filldraw[fill=black](13*0.8,0.7*-30)circle(0.1);        \draw(11*0.8,0.7*-26)--(5*0.8,0.7*-28);        \draw(11*0.8,0.7*-26)--(17*0.8,0.7*-28);        \draw(11*0.8,0.7*-26)--(11*0.8,0.7*-28);        \draw(11*0.8,0.7*-28)--(9*0.8,0.7*-30);        \draw(11*0.8,0.7*-28)--(13*0.8,0.7*-30);        \filldraw[fill=red](14*0.8,0.7*-30)circle(0.1);        \filldraw[fill=red](14*0.8,0.7*-32)circle(0.1);        \filldraw[fill=red](13*0.8,0.7*-32)circle(0.1);        \filldraw[fill=red](15*0.8,0.7*-32)circle(0.1);        \filldraw[fill=red](14*0.8,0.7*-34)circle(0.1);        \filldraw[fill=red](20*0.8,0.7*-30)circle(0.1);        \filldraw[fill=red](17*0.8,0.7*-30)circle(0.1);        \filldraw[fill=red](17*0.8,0.7*-32)circle(0.1);        \filldraw[fill=red](16*0.8,0.7*-32)circle(0.1);        \filldraw[fill=red](18*0.8,0.7*-32)circle(0.1);        \filldraw[fill=red](17*0.8,0.7*-34)circle(0.1);        \filldraw[fill=red](19*0.8,0.7*-34)circle(0.1);        \filldraw[fill=red](15*0.8,0.7*-36)circle(0.1);        \filldraw[fill=red](17*0.8,0.7*-36)circle(0.1);        \filldraw[fill=red](19*0.8,0.7*-36)circle(0.1);        \filldraw[fill=red](19*0.8,0.7*-38)circle(0.1);        \draw[red] (17*0.8,0.7*-28)--(17*0.8,0.7*-30);        \draw[red] (17*0.8,0.7*-28)--(14*0.8,0.7*-30);        \draw[red] (17*0.8,0.7*-28)--(20*0.8,0.7*-30);        \draw[red] (17*0.8,0.7*-30)--(16*0.8,0.7*-32);        \draw[red] (17*0.8,0.7*-30)--(18*0.8,0.7*-32);        \draw[red] (17*0.8,0.7*-30)--(17*0.8,0.7*-32);        \draw[red] (14*0.8,0.7*-30)--(15*0.8,0.7*-32);        \draw[red] (14*0.8,0.7*-30)--(14*0.8,0.7*-32);        \draw[red] (14*0.8,0.7*-30)--(13*0.8,0.7*-32);        \draw[red] (14*0.8,0.7*-30)--(15*0.8,0.7*-32);        \draw[red] (14*0.8,0.7*-32)--(14*0.8,0.7*-34);        \draw[red] (17*0.8,0.7*-34)--(18*0.8,0.7*-32);        \draw[red] (19*0.8,0.7*-34)--(18*0.8,0.7*-32);        \draw[red] (17*0.8,0.7*-34)--(17*0.8,0.7*-36);        \draw[red] (17*0.8,0.7*-34)--(15*0.8,0.7*-36);        \draw[red] (17*0.8,0.7*-34)--(19*0.8,0.7*-36);        \draw[red] (19*0.8,0.7*-36)--(19*0.8,0.7*-38);        \coordinate [label=above:$1$] (x) at (11*0.8,0.7*-26);        \coordinate [label=below:$9$] (x) at (9*0.8,0.7*-30);        \coordinate [label=below:$10$] (x) at (12.4*0.8,0.7*-30);        \coordinate [label={[text=red]above:$2$}] (x) at (14*0.8,0.7*-30);        \coordinate [label={[text=red]right:$3$}] (x) at (17*0.8,0.7*-30);        \coordinate [label={[text=red]below:$8$}] (x) at (20*0.8,0.7*-30);        \coordinate [label={[text=red]below:$6$}] (x) at (19*0.8,0.7*-34);        \coordinate [label={[text=red]left:$4$}] (x) at (17*0.8,0.7*-34);        \coordinate [label={[text=red]below:$7$}] (x) at (19*0.8,0.7*-38);        \coordinate [label={[text=red]below:$5$}] (x) at (14*0.8,0.7*-34);        \coordinate [label=above:$ F$] (x) at (11*0.8,0.7*-6);        \coordinate [label=above:$ F^{(1)}$] (x) at (11*0.8,0.7*-22);        \coordinate [label=above:$ F^{(2)}$] (x) at (11*0.8,0.7*-41);        \coordinate [label=right:$ \Psi_3(F)$] (x) at (11*0.8,0.7*-7);        \coordinate [label=right:$ \Psi_1(F^{(1)})$] (x) at (11*0.8,0.7*-23);                \end{tikzpicture}    \end{center}
  	\caption{ An example of the map $\Gamma$.}\label{FGamma}
  \end{figure}
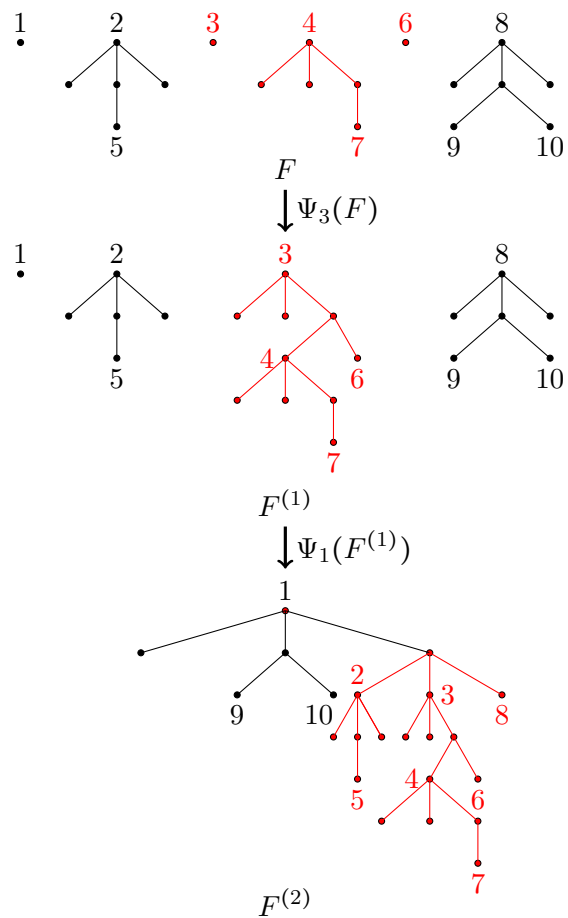

  For example, let $F\in  \mathcal{F}_{10}(3) $ displayed  in Fig.~\ref{FGamma} and let $S=\{1,3\}$. Then we have 
  $\Gamma(F, S) =F^{(2)}$ as shown in Fig.~\ref{FGamma}. 
  
    The following crucial  property of the map $ \Gamma$ follows directly from Lemma \ref{lemPsi}.
  \begin{lemma}\label{propGamma}
  	For any $(F,S)\in \mathcal{Y}_{n}(k) $,   we have 	
  	$F'=\Gamma(F,S)\in  \mathcal{F}_{n}(k) $ with $$\lleaf(F')-\si(F')=	\lleaf(F)-\si(F)+|S|.$$
  	Moreover,  we have $F\in \overline{\mathcal{F}}_n(k)$  if and only if $F'\in  \overline{\mathcal{F}}_n(k)$. 
  	
  \end{lemma}
  
  \noindent{\bf The transformation $\beta$.}\\
  Given a pair $(F,S)  $ where $F=(T_1, T_2, \ldots, T_m)\in \mathcal{F}_n(k)$ with $\rleaf(F)>0$ and $S\subseteq \Si^*(F)$, we define the  transformation $\beta(F,S)$  as follows.  Find the \red{\bf smallest} value $x$ of the set of the labels of the removable old leaves and the removable young leaves. Suppose  that the node with label $x$ is in $ T_i$. Let $F'=\Psi_x(F)$ and $S'=\{y\}\cup S$, where $y$ is the label of the root of $T_i$.  Set $\beta(F,S)=(F', S')$. 
  
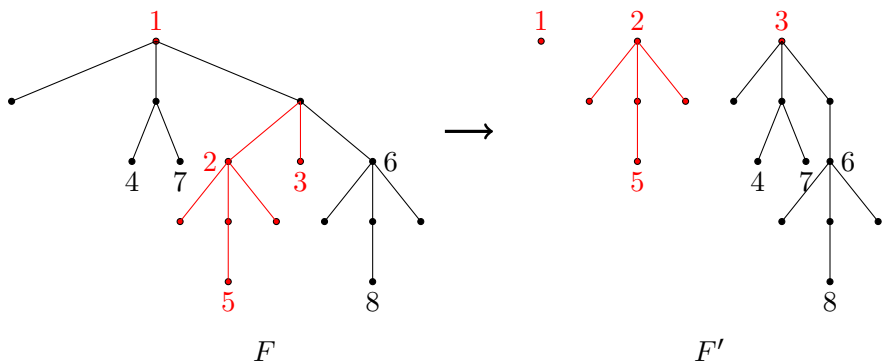
\begin{figure}    \begin{center}	\begin{tikzpicture}[font =\small , scale = 0.4]        \filldraw[fill=red](6*0.80,0)circle(0.1);        \filldraw[fill=black](0*0.80,-2)circle(0.1);        \filldraw[fill=black](6*0.80,-2)circle(0.1);        \filldraw[fill=black](12*0.80,-2)circle(0.1);        \filldraw[fill=black](7*0.80,-4)circle(0.1);        \filldraw[fill=black](5*0.80,-4)circle(0.1);        \filldraw[fill=red](9*0.80,-4)circle(0.1);        \filldraw[fill=red](12*0.80,-4)circle(0.1);        \filldraw[fill=black](15*0.80,-4)circle(0.1);        \filldraw[fill=red](7*0.80,-6)circle(0.1);        \filldraw[fill=red](9*0.80,-6)circle(0.1);        \filldraw[fill=red](11*0.80,-6)circle(0.1);        \filldraw[fill=red](9*0.80,-8)circle(0.1);        \filldraw[fill=black](13*0.80,-6)circle(0.1);        \filldraw[fill=black](15*0.80,-6)circle(0.1);        \filldraw[fill=black](17*0.80,-6)circle(0.1);        \filldraw[fill=black](15*0.80,-8)circle(0.1);        \draw (6*0.80,0)--(0*0.80,-2);        \draw (6*0.80,0)--(6*0.80,-2);        \draw (6*0.80,0)--(12*0.80,-2);        \draw (5*0.80,-4)--(6*0.80,-2);        \draw (7*0.80,-4)--(6*0.80,-2);        \draw (15*0.80,-4)--(12*0.80,-2);        \draw (15*0.80,-4)--(13*0.80,-6);        \draw (15*0.80,-4)--(17*0.80,-6);        \draw (15*0.80,-4)--(15*0.80,-6);        \draw (15*0.80,-8)--(15*0.80,-6);        \draw[red] (12*0.80,-2)--(9*0.80,-4);        \draw[red] (12*0.80,-2)--(12*0.80,-4);        \draw[red] (9*0.80,-4)--(7*0.80,-6);        \draw[red] (9*0.80,-4)--(9*0.80,-6);        \draw[red] (9*0.80,-4)--(11*0.80,-6);        \draw[red] (9*0.80,-8)--(9*0.80,-6);        \coordinate [label={[text=red]above:$1$}] (x) at (6*0.80,0);         \coordinate [label={[text=red]below:$3$}] (x) at (12*0.80,-4);         \coordinate [label={[text=red]left:$2$}] (x) at (9*0.80,-4);         \coordinate [label={[text=red]below:$5$}] (x) at (9*0.80,-8);        \coordinate [label={[text=black]below:$4$}] (x) at (5*0.80,-4);        \coordinate [label={[text=black]below:$7$}] (x) at (7*0.80,-4);        \coordinate [label={[text=black]below:$8$}] (x) at (15*0.80,-8);        \coordinate [label={[text=black]right:$6$}] (x) at (15*0.80,-4);         \draw[line width=1.2pt][->] (18*0.80,-3) -- (20*0.80,-3);                \filldraw[fill=red](22*0.80,0)circle(0.1);        \filldraw[fill=red](26*0.80,0)circle(0.1);        \filldraw[fill=red](24*0.80,-2)circle(0.1);        \filldraw[fill=red](26*0.80,-4)circle(0.1);        \filldraw[fill=red](26*0.80,-2)circle(0.1);        \filldraw[fill=red](28*0.80,-2)circle(0.1);        \draw[red] (26*0.80,0)--(26*0.80,-2);        \draw[red] (26*0.80,-2)--(26*0.80,-4);        \draw[red] (26*0.80,0)--(24*0.80,-2);        \draw[red] (26*0.80,0)--(28*0.80,-2);        \coordinate [label={[text=red]above:$1$}] (x) at (22*0.80,0);        \coordinate [label={[text=red]above:$2$}] (x) at (26*0.80,0);        \coordinate [label={[text=red]below:$5$}] (x) at (26*0.80,-4);        \filldraw[fill=red](32*0.80,0)circle(0.1);        \coordinate [label={[text=red]above:$3$}] (x) at (32*0.80,0);                \filldraw[fill=black](30*0.80,-2)circle(0.1);        \filldraw[fill=black](32*0.80,-2)circle(0.1);        \filldraw[fill=black](34*0.80,-2)circle(0.1);        \filldraw[fill=black](31*0.80,-4)circle(0.1);        \filldraw[fill=black](33*0.80,-4)circle(0.1);        \filldraw[fill=black](34*0.80,-4)circle(0.1);        \filldraw[fill=black](32*0.80,-6)circle(0.1);        \filldraw[fill=black](36*0.80,-6)circle(0.1);        \filldraw[fill=black](34*0.80,-6)circle(0.1);        \filldraw[fill=black](34*0.80,-8)circle(0.1);        \draw (32*0.80,0)--(30*0.80,-2);        \draw (32*0.80,0)--(32*0.80,-2);        \draw (32*0.80,0)--(34*0.80,-2);        \draw (32*0.80,-2)--(31*0.80,-4);        \draw (32*0.80,-2)--(33*0.80,-4);        \draw (34*0.80,-4)--(34*0.80,-2);        \draw (32*0.80,-6)--(34*0.80,-4);        \draw (34*0.80,-6)--(34*0.80,-4);        \draw (36*0.80,-6)--(34*0.80,-4);        \draw (34*0.80,-8)--(34*0.80,-6);        \coordinate [label={[text=black]below:$4$}] (x) at (31*0.80,-4);        \coordinate [label={[text=black]below:$7$}] (x) at (33*0.80,-4);        \coordinate [label={[text=black]below:$8$}] (x) at (34*0.80,-8);        \coordinate [label={[text=black]right:$6$}] (x) at (34*0.80,-4);        \coordinate [label=above:$ F$] (x) at (10.5*0.80,-11);        \coordinate [label=above:$ F'$] (x) at (29*0.80,-11);                \end{tikzpicture}    \end{center}
 \caption{An example of the transformation $\beta$.}\label{beta}
\end{figure}

  For example, let $F$ be the forest on the left of Fig.~\ref{beta} and let $S=\emptyset$.  Then we will get the pair $\beta(F,S)=(F', S')$, where $F'=\Psi_3(F)$ is the forest on the right of Fig.~\ref{beta} and $S'=\{1\}$.

  Let $i$ be a nonnegative integer. Denote by 
  $\beta^i(F, S)$ the resulting pair  obtained after the $i$-th application of $\beta$ to  the pair $(F,S)$. Clearly, we have $\beta^0(F,S)=(F,S)$. 
    
 \noindent{\bf The map $\Gamma'$ from $ \mathcal{F}_{n}(k)$ to $  \mathcal{Y}_{n}(k)$. }\\
  Let $F\in \mathcal{F}_{n}(k) $.  Let $(F^{(0)}, S^{(0)})=(F, \emptyset)$.   
 Relation~\eqref{eqPsi3}
 ensures that after applying finitely many iterations of $\beta$ to $(F^{(0)},S^{(0)})$, 
  we will reach a pair $(F^{(t)},S^{(t)})$ with $\rleaf(F^{(t)})=0$. Set $\Gamma'(F)=(F^{(t)},S^{(t)})$.

 From Lemma~\ref{lemPsi} and  the definition of the  transformation $\Psi_x$, it follows that $(F^{(t)}, S^{(t)})\in\mathcal{Y}_{n}(k)$ and hence the map $\Gamma'$ is well-defined.

   Regarding the transformation $\beta$, we have the following observation, which can be checked routinely.
  
  \begin{observation}\label{obser:beta}
  Let $F\in \mathcal{F}_n(k)$ with $\rleaf(F)>0$. Let $\beta^{i-1}(F,\emptyset)=(F', S')$  with  $\rleaf(F')>0$ for some $i\geq1$.  Suppose that  $\beta^{i}(F,\emptyset)=(F'', S'')=(\Psi_x(F'), S'\cup \{y\})$.    Then   there does not exist  any removable leaves  in $F''$ that are located to the left of the singleton with label $y$.
  \end{observation}

\begin{proposition}\label{Gamma0}
	For $n,k\geq 1$, the maps $\Gamma$ and $\Gamma'$ induce a bijection between $\overline{\mathcal{Y}}_{n}(k)$ and $\overline{\mathcal{F}}_n(k)$ such that for any $(F,S)\in \overline{\mathcal{Y}}_n(k)$, we   have 	
	$F'=\Gamma(F,S)\in \overline{\mathcal{F}}_{n}(k) $ with $$\lleaf(F')-\si(F')=	\lleaf(F)-\si(F)+|S|.$$
\end{proposition}
\pf In view of Lemma~\ref{propGamma}, it suffices  to show that the maps $\Gamma$ and $\Gamma'$ are inverses of each other.   First, we proceed to show that $\Gamma'
$ is the inverse 
of the map  $\Gamma$, that is, $\Gamma'(\Gamma(F,S))=(F,S)$ for all $(F,S)\in \overline{\mathcal{Y}}_n(k)$. To this end, it suffices to show that $\beta(\alpha^{j}(F,S))=\alpha^{j-1}(F,S)$.
Let $(F', S')=\alpha^{j-1}(F,S)$ with $F'=(T_1, T_2, \ldots, T_m)$. 
Suppose that at the $j$-th application of $\alpha$ to $(F,S)$, $x$ is  the greatest value of the set $S'$ and $T_i$ is the singleton with label $x$.
Suppose  $s>i$ is the least integer  such that   $T_s$ is either a singleton  or $T_s=T_m$. Assume that  the root of $T_s$ is labeled by $y$. Then  
by the definition of the transformation $\Psi_x$, one can easily check that in the forest $\Psi_x(F')$,  the tree rooted at the node with label $x$ contains a removable old/young leaf  with label   $y$ such that  $y$ is the smallest value of the set of the labels of removable old leaves and removable young leaves of 
$\Psi_x(F')$.  Again by the definition of $\Psi_x$ and $\Psi_y$, it is routine to check that  $\Psi_y(\Psi_x(F'))=F'$. 
Therefore, we deduce that
$$\beta(\alpha^{j}(F,S))=\beta(\alpha(F',S'))=\beta(\Psi_x(F'), S'-\{x\} )=(\Psi_y(\Psi_x(F')), S') =(F',S'),$$ as desired. So far, we have concluded that  
$\Gamma'
$ is the inverse 
of the map  $\Gamma$.

Next we aim to show that $\Gamma$ is the inverse 
of the map  $\Gamma'$, that is, $\Gamma(\Gamma'(F))=F$ for all $F\in \overline{\mathcal{F}}_n(k)$. To this end, it suffices to show that $\alpha(\beta^{j}(F,\emptyset))=\beta^{j-1}(F,\emptyset)$.
Let $(F', S')=\beta^{j-1}(F,\emptyset)$. 
Suppose that at the $j$-th application of $\beta$ to $(F, \emptyset)$,   the selected removable leaf has label $x$. If the leaf $x$ is in  $T_i$ whose root has label $y$.  Then  
by the definition of the transformation $\Psi_y$, one can easily check that   the forest $\Psi_y(F')$ contains a singleton    with label $y$.  Recall that $\beta(F', S')=(\Psi_x(F'), S'\cup \{y\})$. By Observation~\ref{obser:beta}, $y$ is the greatest value in $S'\cup \{y\}$.  Then we have $\alpha(\Psi_x(F'), S'\cup \{y\})=(\Psi_y(\Psi_x(F')), S')$. Again by the definition of $\Psi_y$ and $\Psi_x$, it is routine to check that  $\Psi_y(\Psi_x(F'))=F'$. 
Therefore, we deduce that
$$\alpha(\beta^{j}(F,\emptyset))=\alpha(\beta(F',S'))= (\Psi_y(\Psi_x(F')), S') =(F',S'),$$ as desired. Thus, we conclude  that  
$\Gamma
$ is the inverse 
of the map  $\Gamma'$, completing the proof. \qed 

The following proposition can be verified by the same reasoning as in the proof of Proposition~\ref{Gamma0} and the proof is omitted here. 
\begin{proposition}\label{Gamma}
	For $n,k\geq 1$, the maps $\Gamma$ and $\Gamma'$ induce a bijection between $\widehat{\mathcal{Y}}_{n}(k)$ and $\widehat{\mathcal{F}}_n(k)$ such that for any $(F,S)\in \widehat{\mathcal{Y}}_n(k)$, we   have 	
	$F'=\Gamma(F,S)\in \widehat{\mathcal{F}}_{n}(k) $ with $$\lleaf(F')-\si(F')=	\lleaf(F)-\si(F)+|S|. $$
\end{proposition}

Combining  Propositions \ref{Theta0}, \ref{Theta},    \ref{Gamma0} and \ref{Gamma}, we are  led to the following results. 
\begin{proposition}\label{Mainbijection0}
	For $n,k\geq 1$, the map  $\Gamma\circ \Theta$ induces a one-to-one correspondence  between $\overline{\mathcal{X}}_{n}(k)$ and $\overline{\mathcal{F}}_n(k)$ such that for any $(F,S)\in \overline{\mathcal{X}}_n(k)$, we   have 	
	$F'=\Gamma( \Theta(F,S))\in \overline{\mathcal{F}}_{n}(k) $ with $$\lleaf(F')-\si(F')=	\lleaf(F)-\si(F)+|S|.$$
\end{proposition}

\begin{proposition}\label{Mainbijection}
	For $n,k\geq 1$, the map  $\Gamma\circ \Theta$ induces a one-to-one correspondence  between $\widehat{\mathcal{X}}_{n}(k)$ and $\widehat{\mathcal{F}}_n(k)$ such that for any $(F,S)\in \widehat{\mathcal{X}}_n(k)$, we   have 	
	$F'=\Gamma(\Theta(F,S))\in \widehat{\mathcal{F}}_{n}(k) $ with $$\lleaf(F')-\si(F')=	\lleaf(F)-\si(F)+|S|.$$
\end{proposition}

   \subsection{Finishing the proof of Theorem~\ref{gammaa}}
   
    Now we are ready to complete the proof of Theorem~\ref{gammaa}.  In view of~\eqref{eq:ab}, it remains to show that  
    \begin{equation}\label{eqa2}
 	  	\sum\limits_{F\in \overline{\mathcal{F}}_{n}(k)}x^{\lleaf(F)-\si(F)}=\sum_{i=0}^{\lfloor{n-1\over 2}\rfloor} \overline{\gamma}_{n, k,i}x^i(1+x)^{n-1-2i}
 	  \end{equation}
	  and
	   \begin{equation}\label{eqb2}
\sum\limits_{F\in \widehat{\mathcal{F}}_{n}(k)}x^{\lleaf(F)-\si(F)}=\sum_{i=1}^{\lfloor{n\over 2}\rfloor} \widehat{\gamma}_{n,k, i}x^i(1+x)^{n-2i}.
  	\end{equation}
 
   First we prove~\eqref{eqa2}. 
   By Lemma {\ref{relation1}}, we have 
   \begin{equation}\label{relation20}
   	\oint(F)+\si(F)=n-2\oleaf(F)  	
   \end{equation}
   for any $F\in  \overline{\mathcal{F}}^*_{n}(k)$. 
   Note that for any $F\in \mathcal{F}_n(k)$, we have
   \begin{equation}\label{relation:oys}
   \oleaf(F)+\yleaf(F)+\si(F)=\lleaf(F).
   \end{equation}
   Hence,   for any $F\in  \overline{\mathcal{F}}^*_{n}(k)$, we have   \begin{equation}\label{relation30} 
   	\oleaf(F)+\si(F)=\lleaf(F). 
   \end{equation}
 For any  $F\in \overline{\mathcal{F}}_n(k)$,  one can easily check that $\oint^*(F)+\si^*(F)= \oint(F)+\si(F)-1$, where $\oint^*(F)=|\Oint^*(F)|$ and  $\si^*(F)=|\Si^*(F)|$. This combined with~\eqref{relation20} yields that
   \begin{equation}\label{relation4}
  	\oint^*(F)+ \si^*(F)=n-1-2\oleaf(F)  	
  \end{equation}
for any $F\in \overline{\mathcal{F}}^*_{n}(k)$.
   
   Note that
   $\overline{\gamma}_{n, k,i}$ is  the number of forsets in  $\overline{\mathcal{F}}^*_n(k)$  with $i$ old leaves.
   In view of Proposition \ref{Mainbijection0}, we  deduce that
   $$
   \begin{array}{lll}
   	\sum\limits_{F\in \overline{\mathcal{F}}_n(k)}x^{\lleaf(F)-\si(F)}&=&  \sum\limits_{(F', S')\in \overline{\mathcal{X}}_{n}(k)}x^{\lleaf(F')-\si(F')+|S'|}\\
   	&=&  \sum\limits_{F'\in  \overline{\mathcal{F}}^*_{n}(k)}\,\,\sum\limits_{S'\subseteq \Oint^*(F')\bigcup \Si^*(F')}x^{\lleaf(F')-\si(F')+|S'|}\\
   	&=&  \sum\limits_{F'\in  \overline{\mathcal{F}}^*_{n}(k)}\,\,\sum\limits_{S'\subseteq \Oint^*(F')\bigcup \Si^*(F')}x^{\oleaf(F')+|S'|} \\
       	&=&  \sum\limits_{F'\in  \overline{\mathcal{F}}^*_{n}(k)} x^{\oleaf(F') }(1+x)^{\oint^*(F')+\si^*(F')} \\
   	&=&  \sum\limits_{F'\in  \overline{\mathcal{F}}^*_{n}(k)} x^{\oleaf(F')} (1+x)^{n-1-2\oleaf(F')}\\
   	&=&  \sum\limits_{i=0}^{\lfloor{n-1\over 2}\rfloor} \overline{\gamma}_{n, k,i}x^i(1+x)^{n-1-2i}, 
   \end{array}
   $$   
   where the third equality follows from~\eqref{relation30} and the fifth equality follows from 
   \eqref{relation4}.  This completes the proof of~\eqref{eqa2}.

Next we prove~\eqref{eqb2}.   
By Lemma {\ref{relation1}}, we have 
 \begin{equation}\label{relation2}
 \oint(F)+\si(F)=n-2\oleaf(F)  	
 	\end{equation}
 for any $F\in  \widehat{\mathcal{F}}^*_{n}(k)$.
 It follows from~\eqref{relation:oys} that 
  \begin{equation}\label{relation3} 
 	\oleaf(F)+\si(F)=\lleaf(F)
 	\end{equation}
	for any $F\in  \widehat{\mathcal{F}}^*_{n}(k)$. 
  Note that for any $F\in \widehat{\mathcal{F}}_{n}(k)$, we have $\Si^*(F)=\Si(F)$. 
  
  It is clear that
$\widehat{\gamma}_{n,k, i}$ is  the number of forests in  $\widehat{\mathcal{F}}^*_n(k)$ with $i$ old   leaves.
In view of Proposition \ref{Mainbijection}, we  derive that
\begin{align*}
\sum\limits_{F\in \widehat{\mathcal{F}}_n(k)}x^{\lleaf(F)-\si(F)}&=  \sum\limits_{(F', S')\in \widehat{\mathcal{X}}_{n}(k)}x^{\lleaf(F')-\si(F')+|S'|}\\
&=  \sum\limits_{F'\in  \widehat{\mathcal{F}}^*_{n}(k)}\,\,\sum\limits_{S'\subseteq \Oint(F')\bigcup \Si^*(F')}x^{\lleaf(F')-\si(F')+|S'|}\\
&=  \sum\limits_{F'\in  \widehat{\mathcal{F}}^*_{n}(k)}\,\,\sum\limits_{S'\subseteq \Oint(F')\bigcup \Si^*(F')}x^{\oleaf(F')+|S'|} \\
&=  \sum\limits_{F'\in  \widehat{\mathcal{F}}^*_{n}(k)}\,\,\sum\limits_{S'\subseteq \Oint(F')\bigcup \Si(F')}x^{\oleaf(F)+|S'|} \\
&=  \sum\limits_{F'\in  \widehat{\mathcal{F}}^*_{n}(k)} x^{\oleaf(F') }(1+x)^{\oint(F')+\si(F')} \\
&=  \sum\limits_{F'\in  \widehat{\mathcal{F}}^*_{n}(k)} x^{\oleaf(F')} (1+x)^{n-2\oleaf(F')}\\
&=\sum_{i=1}^{\lfloor{n\over 2}\rfloor} \widehat{\gamma}_{n,k, i}x^i(1+x)^{n-2i},
\end{align*}
 where the third equality follows from  \eqref{relation3} and the sixth  equality follows from~\eqref{relation2}. 
This proves~\eqref{eqb2}, which completes the proof of  Theorem~\ref{gammaa}.

 \section*{Acknowledgments}
 
This work 
was supported by the National  Science Foundation of China (grants 12471318, 12071440, 12322115 \& 12271301) 
 and the Fundamental Research Funds for the Central Universities.


\end{document}